\documentclass[11pt]{amsart}
\usepackage{amsfonts}
\usepackage{amssymb}
\usepackage{graphicx}
\usepackage{amsmath}
\usepackage{latexsym}
\usepackage{amscd}
\usepackage{mathrsfs}
\usepackage{enumitem} 
\usepackage{braket}
\usepackage[margin=1.4in]{geometry}
\usepackage{bm}
\usepackage{color}
\usepackage{marginnote}
\usepackage{verbatim}

\usepackage{imakeidx}
\makeindex[title=Index of Notation]

\usepackage{hyperref}

\allowdisplaybreaks 

\newtheorem{theorem}{Theorem}[section]
\newtheorem{corollary}[theorem]{Corollary}
\newtheorem{lemma}[theorem]{Lemma}
\newtheorem{proposition}[theorem]{Proposition}
\theoremstyle{definition}
\newtheorem{definition}[theorem]{Definition}
\newtheorem{example}[theorem]{Example}
\newtheorem{remark}[theorem]{Remark}
\newtheorem{defnlem}[theorem]{Definition/Lemma}

\newcommand{\vp}{\varphi}
\newcommand{\HH}{\mathcal{H}}
\newcommand{\R}{\mathbb{R}}
\newcommand{\N}{\mathbb{N}}
\newcommand{\Z}{\mathbb{Z}}
\newcommand{\C}{\mathcal{C}}
\newcommand{\SSS}{\mathbb{S}}
\newcommand{\LL}{\mathcal{L}}
\newcommand{\IV}{\mathcal{IV}}
\newcommand{\loc}{\mathrm{loc}}
\newcommand{\fint}{\mathrel{\int\!\!\!\!\!\!-}}
\newcommand{\ssubset}{\subset\subset}

\DeclareMathOperator{\spt}{spt}
\DeclareMathOperator{\sing}{sing}
\DeclareMathOperator{\proj}{Proj}
\DeclareMathOperator{\reg}{reg}
\DeclareMathOperator{\id}{id}
\DeclareMathOperator{\indec}{Indec}

\DeclareMathOperator{\lip}{Lip}
\DeclareMathOperator{\Div}{div}
\DeclareMathOperator{\graph}{graph}
\DeclareMathOperator{\esssup}{ess\ sup}
\begin{document}

\begin{abstract}
    We establish a local Harnack inequality in a neighborhood of an indecomposable singular point of a stationary integral varifold. Extending the method of Gr\"uter and Widman \cite{gruter1982green}, we construct the Green function on a stationary integral varifold with Euclidean volume growth, allowing the pole to be any point in the support. Using the local Harnack inequality, we show that if a sequence of stationary varifolds converges with multiplicity one, then the corresponding Green functions converge as well. As applications, we determine the asymptotic behavior of the Green function both near an indecomposable pole and at infinity. We further establish global lower and upper bounds for the Green function and present an application of these estimates. Finally, we analyze the behavior of the Green function when the varifold is decomposable at a point, including examples illustrating the possible phenomena.
\end{abstract}
\title{Green functions on stationary varifolds}
\author{Yifan Guo}
\address{Courant Institute, New York University, 251 Mercer St, New York, NY, USA, 10012}
\email{yg3991@nyu.edu}
\maketitle
\section{Introduction}
There has been a long history of studying Green functions on complete noncompact manifolds, particularly on those with nonnegative Ricci curvature (cf. \cite{CLY,CM1,D,littman1963regular,LT,LY,Var}). In this work, we investigate Green functions on singular minimal submanifolds, or more precisely, on stationary varifolds. Our focus is on their convergence, asymptotic behavior near the pole and at infinity, as well as global upper and lower bounds.

Many classical results rely on the Bochner formula, the maximum principle, and lower Ricci curvature bounds—tools that are not available in the setting of stationary varifolds. Our approach is closer in spirit to that of \cite{D}, where convergence of Green functions is established on Ricci limit spaces. A key ingredient in proving our convergence theorem is a local Harnack inequality, which we establish near an indecomposable point of a stationary varifold.

\subsection{Indecomposability and local Harnack inequality}
We consider integral $n$-varifolds $V \in \IV_n(U)$ that are stationary in an open set $U \subset \mathbb{R}^{n+k}$. 
For $x \in \spt V$, we denote by $\tan(V,x)$ the set of all tangent cones of $V$ at $x$. 
If $V$ has Euclidean volume growth, we write $\tan(V,\infty)$ for the set of all tangent cones at infinity.

Indecomposability plays an essential role in the validity of a Harnack inequality. Indeed, if $V$ is the union of two stationary varifolds $V_1$ and $V_2$, then the function defined by $u=1$ on $V_1$ and $u=0$ on $V_2$ is harmonic on $V$, yet no Harnack inequality can hold. 
Motivated by this observation, we introduce the following notions of indecomposability for stationary varifolds and for points.

\begin{definition}\label{definition.indecomposable}
    Let $V\in \IV_n(U)$ be stationary. We say that $V$ is \textbf{indecomposable} in $U$ if for any $V_1,V_2\in \IV_n(U)$ such that $V=V_1+V_2$ and $\delta V_1=\delta V_2=0$, we have either $V_1=0$ or $V_2=0$. 
    
    Let $x\in \spt V$. We say that $V$ is \textbf{indecomposable  at }$x$ if for any $\C\in \tan(V,x)$, $\C$ is indecomposable. We write \[\indec V=\{x\in \spt V: x\text{ is indecomposable}\}.\]

    If $V$ has Euclidean volume growth, we say that $V$ is \textbf{indecomposable  at }$\infty$ if any $\C\in \tan(V,\infty)$ is indecomposable.
\end{definition}
Related notions of indecomposability have been considered in \cite[Definition~2.15]{Modino} and \cite[Definition~6.2]{Me2weaklydifferentiable}, although our definitions differ slightly from those references; see Remark~\ref{different.notion.indecomposable}.

Under mild assumptions on $\sing V$, indecomposability is equivalent to multiplicity one and connectivity of the regular part.

\begin{lemma}   
Let $V \in \IV_n(U)$ be a stationary varifold.  

(i) If $\mathcal{H}^n(\sing V) = 0$, $V$ has multiplicity one and $\reg V$ is connected, then $V$ is indecomposable.  

(ii) Conversely, if $\mathcal{H}^{n-1}(\sing V) = 0$ and $V$ is indecomposable, then $V$ has multiplicity one and $\reg V$ is connected.

\end{lemma}
The following example shows that the assumption $\HH^{n-1}(\sing V)=0$ in (ii) cannot be relaxed.

\begin{example}\label{triple.junction.intro}
    Let $H_1,H_2,H_3$ be three half $n$-planes in $\R^{n+k}$ with $\partial H_1=\partial H_2=\partial H_3=\R^{n-1}$ and $H_i$ and $H_j$ have angle $120^\circ$ for $i\neq j$. Let $V=H_1+H_2+H_3$. Then $\sing V=\partial H_1$ and $V$ is indecomposable by Lemma \ref{triple.junction.indecomposable} but $\reg V$ consists of 3 components.
\end{example}

The first main result of this article is the following local relative isoperimetric inequality and local Harnack inequality near an indecomposable point. In \cite[Theorem~2, Theorem~5]{bombieri1972harnack}, Bombieri and Giusti establish the relative isoperimetric inequality and the Harnack inequality for area-minimizing boundaries. In their setting no indecomposability assumption is required, since area-minimizing boundaries are automatically indecomposable \cite[Theorem~1]{bombieri1972harnack}. 

Although our proof is largely similar to the argument in \cite{bombieri1972harnack}, a key new feature is that if $V$ is indecomposable at $x$, then the relative isoperimetric inequality and the Harnack inequality hold in a sufficiently small ball for \emph{all nearby} stationary varifolds.

To state the theorem, we recall the notion of the distributional boundary of a set $E$ in $V$, introduced in \cite[Definition~5.1]{Me2weaklydifferentiable}:
\[
V\partial E=(\delta V)\llcorner E-\delta (V\llcorner E).
\]

\begin{theorem}
	Let $V\in \IV_n(U)$ be stationary and $x\in \indec V$. Then there exists $\beta_0=\beta_0(n,V,x)\in(0,1)$, $r_0=r_0(n,V,x)>0$, $\epsilon_0=\epsilon_0(n,V,x)\in(0,1)$ with the following significance.  For any $r<r_0$, $V'\in \IV_n(B_r(x))$ stationary  with $\mathbf{F}_{B_r(x)}(V,V')<\epsilon_0 r^n$ and $\spt V'\cap B_{\beta_0 r/2}(x)\neq\emptyset$ we have the following. 
    
    (i) If $E\subset \R^{n+k}$ is $\|V'\|$-measurable, then
	\[\|V'\partial E\|(B_r(x))=\|V'\partial E^c\|(B_r(x))\ge \frac{1}{2\gamma}\min\{\|V'\|( E\cap B_{\beta_0 r}(x)),\|V'\|(E^c\cap B_{\beta_0 r}(x))\}^{\frac{n-1}{n}}\]
    where $\gamma=\gamma(n,k)>0$.
    
    (ii) If $0\le u\in W^{1,2}(B_{r}(x),V')$ satisfies $-\Delta_{V'}u=f$ in $B_{r}(x)$ for $f\in L^{q}(B_r(x),\|V'\|)$ with $q>n/2$, then we have 
    \[
    \sup_{\spt V'\cap B_{\beta_0 r}(x)}u\le C(\inf_{\spt V'\cap B_{\beta_0 r}(x)} u+r^{2-\frac{n}{q}}\|f\|_{L^{q}(B_r(x),\|V'\|)})
    \]
    for some $C=C(n,k,q,\Lambda,\beta_0)>0$ where $\Lambda=\omega_n^{-1}r_0^{-n}\|V'\|(B_{r_0}(x))$ where $\omega_n=\HH^{n}(B_1)$.
\end{theorem}

\subsection{The Green function on a stationary varifold}Let $n \ge 3$ and $k \ge 1$. We construct the Green function on a stationary integral varifold with Euclidean volume growth by adapting the method of Gr\"uter and Widman \cite{gruter1982green}. In their work, they construct the Green function for uniformly elliptic operators with measurable coefficients on Euclidean domains, relying on analytic tools such as the Sobolev inequality and the mean value inequality for positive harmonic functions. These tools are also available on stationary integral varifolds with Euclidean volume growth, allowing us to extend their approach to this more general setting. Among other results, we obtain the following existence theorem.
\begin{theorem}
	Let $V\in \IV_n(\R^{n+k})$ be stationary with Euclidean volume growth. There exists a function $G:\spt V\times \spt V\to [0,\infty]$
		satisfying: for any $x\in \spt V$, $r>0$ we have 
		\begin{align*}
			&G(x,\cdot)\in W_\loc^{1,2}( \R^{n+k}\setminus B_{r}(x),V)\cap W_\loc^{1,p}(V)\quad \forall p\in [1,\frac{n}{n-1})\\
			&G(x,y)\to0\quad\text{as }|y|\to\infty
		\end{align*}
		and  $\Delta_y G(x,y)=-\delta_x(y)$ in the sense that for any $\phi\in C_c^\infty(\R^{n+k})$, 
		\[
			\int\nabla _y^V G(x,y) \cdot \nabla^V\phi(y)d\|V\|(y)=\phi(x).
		\]
\end{theorem}There is also a version of Dirichlet Green function $G_\Omega(x,y)$ of $V$ restricted to a bounded open set $\Omega\subset \R^{n+k}$. See Proposition~\ref{DiriGreen} for more details.
We note that the Green function is defined with the pole $x$ be any point in $\spt V$. However, some properties may not be true everywhere e.g. we only have $G(x,y)=G(y,x)$ for $\|V\|$-a.e. $x,y$ or any $x,y\in \reg V\cup \indec V$. See Theorem~\ref{globalGreenvarifolds} for more details and Example~\ref{2planes} for a Green function which is not symmetric at all points $x,y\in \spt V$. 

An important example of the Green function is the Green function on a cone.
\begin{example}
	Let $\C\in \IV_n(\R^{n+k})$ be a stationary cone. Then
	\[G(0,y)=\frac{1}{n(n-2)\omega_n\Theta(\C,0)}|y|^{2-n}.\]
\end{example}

\subsection{Convergence of the Green functions and its consequences}
Using the local Harnack inequality, we establish the following convergence of Green functions under convergence of multiplicity one stationary varifolds. This is the main tool to study the Green function on a stationary varifold. 
We  record the following assumptions on stationary integral varifolds $V\in \IV_n(U)$, to be used later: 
\begin{equation}
	\text{for any } W\ssubset U,\quad \HH^{n-2}(\sing V\cap W)<\infty\label{Hn-2(sing)<infty}.
\end{equation}
\begin{theorem}\label{conv}
	Let $V,\{V_i\}_{i=1}^\infty$ be stationary integral $n$-varifolds in $\R^{n+k}$ with Euclidean volume growth. Assume that $V_i \to V$ in $\IV_n(\R^{n+k})$ and that $V$ has multiplicity one and satisfies \eqref{Hn-2(sing)<infty}. Let $\{\Omega_i\}_{i=1}^\infty$ be either an exhaustion of $\R^{n+k}$ by bounded open subsets, or $\Omega_i = \R^{n+k}$ for all $i$. Let $G$ be the Green function of $V$ and $G_i$ be the Green function of $V_i$ on $\Omega_i$. Assume that $\spt V_i\ni x_i\to x_0\in \reg V\cup \indec V$. 
    
    Then we have $G_i(x_i,\cdot)\to G(x_0,\cdot)$ in $L^q$ on compact subsets as $V_i\to V$ for $1\le q<\frac{n}{n-2}$. Moreover, the convergence is uniform in a neighborhood of any $y\in \indec V\setminus \{x_0\}$ and in any small annulus centered at $x_0$.
    
    If $V=\C$ is a cone and $x_0=0$, then $G_i(x_i,\cdot)\to G(0,\cdot)$ uniformly in any annulus centered at $0$.
\end{theorem}
The first consequence of the convergence theorem above is the following.
\begin{theorem}\label{asymp}
    Let $V\in \IV_n(\R^{n+k})$ be stationary with Euclidean volume growth and $\Omega\subset \R^{n+k}$ be a bounded open set. Denote by $G_\Omega(x,y)$ the Dirichlet Green function of $V$ on $\Omega$ and by $G(x,y)$ the Green function of $V$.
    If $x\in \reg V$ or $x\in \indec  V$ such that any $\C\in \tan(V,x)$ satisfies \eqref{Hn-2(sing)<infty}, then we have
    \[\lim_{y\to x}G_\Omega(x,y)|y-x|^{n-2}=\lim_{y\to x}G(x,y)|y-x|^{n-2}=\frac{1}{n(n-2)\omega_n\Theta(V,x)}.
	\]
	If $V$ is indecomposable at $\infty$ and any $\C\in \tan(V,\infty)$ satisfies \eqref{Hn-2(sing)<infty}, then for any $x\in\spt V$
	\[\lim_{|y|\to \infty}G(x,y)|y-x|^{n-2}=\frac{1}{n(n-2)\omega_n\Theta(V,\infty)}.\]
\end{theorem}
With the above theorem, we have the following sharp upper bound of the Green function. 

\begin{corollary}Let $V\in \IV_n(\R^{n+k})$ be stationary with Euclidean volume growth. Denote by $G(x,y)$ the Green function of $V$. Let $x\in \reg V$ or $x\in \indec  V$ such that any $\C\in \tan(V,x)$ satisfies \eqref{Hn-2(sing)<infty}. Then for $\|V\|$-a.e. $y\in \spt V$ or $y\in \reg V\cup \indec V$,  we have
	\begin{equation}\label{upperglobal.intro}
	    G(x,y)\le\frac{1}{n(n-2)\omega_n\Theta(V,x)}|x-y|^{2-n}.
	\end{equation}
	If $\HH^{n-1}(\sing V)=0$ and equality in \eqref{upperglobal.intro} holds at $x\neq y\in \reg V\cup \indec V$, then $V=\eta_{x,1\#}\C_0+V_1$ where $\C_0$ is a stationary cone and $y-x\in \spt \C_0$, $\eta_{x,1}$ is the map $z\mapsto z-x$ and $V_1\in \IV_n(\R^{n+k})$ is stationary with $x\notin \spt V_1$.
\end{corollary}

We remark that the sharp upper bound for the Green function is a well-known result of Cheng, Li and Yau \cite{CLY} on \textit{smooth} minimal submanifolds $M^n\subset \R^{n+k}$. They proved this by comparing the heat kernels on minimal submanifolds with the heat kernel on $\R^n$ and integrating the time variable. One can also directly compare the Green functions on minimal submanifolds with the Green function on $\R^n$ (cf. Lemma \ref{super}). Our contribution lies in extending it to the singular case and discussing the rigidity case when equality holds.

We also have the following global lower bounds of the Green function when every singular point is indecomposable.
\begin{theorem}\label{lower.bound.intro}
    Let $V\in \IV_n(\R^{n+k})$ be stationary with Euclidean volume growth and $x\in \indec V$. Denote by $G(x,y)$ the Green function of $V$. Suppose $\spt V=\reg V\cup\indec V$ and $V$ is indecomposable at $\infty$ and that any $\C\in \tan(V,x)\cup  \tan(V,\infty)$ satisfies \eqref{Hn-2(sing)<infty}. Then there is $c=c(n,V,x)>0$ such that for any $y\in \spt V$
	\[ G(x,y)\ge c(n,V,x)|x-y|^{2-n}.\]
\end{theorem}

We note that the indecomposability assumption in Theorem \ref{conv}, \ref{asymp}, \ref{lower.bound.intro} is necessary as can be seen in Example \ref{2planes}.

As an application, we have the following gradient estimates for entire solutions to the minimal surface equation.
\begin{corollary}
	Let $u:\R^n\to \R$  be an entire solution to the minimal surface equation 
	\[(1+|\nabla u|^2)\Delta u-\sum_{i,j=1}^nu_{x_ix_j}u_{x_i}u_{x_j}=0.\]
	Then there exists $C=C(n,u)>0$ such that for $|x|\ge 1$
	\[|\nabla u(x)|\le C|x|^{n-2}\left(1+\frac{|u(x)-u(0)|}{|x|}\right)^{n-2}.\]
\end{corollary}
\subsection{Green functions on decomposable varifolds}
We present results of Green functions on stationary varifolds with decomposable singular points or decomposable tangent cones at infinity.

The first is an example which is the simplest instance of a decomposable stationary varifold. It exhibits distinct behavior of the Green function depending on whether the intersection of $V_1,V_2$ has codimension one or codimension at least two. In particular, if $V_1$ and $V_2$ intersect along a set of codimension at least two, then the Green function is simply the union of the respective Green functions on $V_1$ and $V_2$ (cf. Corollary \ref{Green.V_1+V_2}). In contrast, if $V_1$ and $V_2$ intersect along a codimension-one set, the Green function propagates from $V_1$ to $V_2$ when the pole lies on $V_1$, and vice versa. 

For $x,y\in \R^{n+k}$, we write the Euclidean Green function on $\R^n$ as
\[\Gamma(x,y)=\frac{1}{n(n-2)\omega_n}|x-y|^{2-n}.\]

\begin{example}\label{2planes}
	Let $V=P_1+P_2$ be the union of two distinct $n$-planes $P_1,P_2$ in $\R^{n+k}$ and $G$ be its Green function. If the pole $x\in P_1\cap P_2$, then for any $y\in P_1\cup P_2$ we have
	\[G(x,y)=\frac{1}{2}\Gamma(x,y).\]
	If $\dim P_1\cap P_2\le n-2$, then
	\[G(x,y)=\left\{\begin{aligned}
		&\Gamma(x,y)&\quad& x,y\in P_1\setminus P_2\text{ or }x,y\in P_2\setminus P_1\\
		&0 &\quad &\text{otherwise.}
	\end{aligned}\right.\] 
    If $\dim P_1\cap P_2=n-1$, then for $x\in P_1^+$
	\[G(x,y)=\left\{\begin{aligned}
		&\Gamma(x,y)-\frac{1}{2} \Gamma (x^*,y)&\quad& y\in P_1^+\\
		&\frac{1}{2}\Gamma(x^*,y)&\quad&y\in P_1^-\cup P_2^+\cup P_2^- .
	\end{aligned}\right.\]
	To explain the notation, since $\dim P_1\cap P_2=n-1$, $(P_1\cup P_2)\setminus(P_1\cap P_2)=P_{1}^+\cup P_1^-\cup P_2^+\cup P_2^-$. We identify $P_i^\pm$ with $\R^n_+=\{(x_1,\dots,x_n)\in \R^n:x_n>0\}$ in the following way. Let $\{e_1\dots,e_{n-1}\}$ be an orthonormal basis of $P_1\cap P_2$ and $\nu_i^\pm$ be the inward conormal of $P_1\cap P_2$ in $P_i^\pm$. We take the basis $\{e_1,\dots,e_{n-1},\nu_i^\pm\}$ on $P_i^\pm$. For any $x=(x_1,\dots,x_n)\in \R^n_+$, we write $x^*=(x_1,\dots,x_{n-1},-x_n)$.
\end{example}
A similar phenomenon also occurs in the presence of triple junction singularities, even though the varifold is indecomposable.
\begin{example}
    Let $V=H_1+H_2+H_3$ be as in Example~\ref{triple.junction.intro} and $G$ be its Green function. 
    For any $x\in H_1\cap H_2\cap H_3$, we have
	\[G(x,y)=\frac{2}{3}\Gamma(x,y).\]
	As in Example \ref{2planes}, we identify $H_i$ with $ \R^n_+$ for $i=1,2,3$.
	If $x\in H_1\setminus (H_1\cap H_2\cap H_3)$ then
	\[G(x,y)=\left\{\begin{aligned}
		&\Gamma(x,y)-\frac{1}{3} \Gamma (x^*,y)&\quad& y\in H_1\\
		&\frac{2}{3}\Gamma(x^*,y)&\quad&\text{on }  y\in H_2\cup  H_3 
	\end{aligned}\right.\]
\end{example}

The following Lemma shows that although indecomposability is sufficient for obtaining the asymptotic behavior of the Green function in Theorem~\ref{asymp}, and is necessary in view of Example \ref{2planes}, there are nevertheless situations in which one can still derive asymptotics even when the tangent cone is decomposable. Examples satisfying the assumptions include the $n$-dimensional catenoid in $\R^{n+1}$ for $n\ge 3$ and the Lawlor's neck introduced in \cite{La}. 
\begin{lemma}
    Let $\Sigma^n\subset \R^{n+k}$ be smooth minimal submanifold with two ends $E_1,E_2$. We assume that for $i=1,2$, $E_i$ is regular in the sense that $E_i$ can be written as a graph of $v_i$ over some $n$-plane $P_i$. See Lemma~\ref{regular.at.infinity.2.planes} and \cite[Definition, page 800]{Sch83} for more details. We also assume that there exist an isometry of $\Sigma$  such that one end is mapped to the other one. Then we know that the tangent cone at infinity of $\Sigma$  is $P_1+ P_2$ which is decomposable but for any  $x\in\Sigma $, we have\[\lim_{|y|\to\infty}|x-y|^{n-2}G(x,y)=\frac{1}{2n(n-2)\omega_n}.\] 
\end{lemma}
The outline of the rest of the paper is as follows. In section \ref{preliminaries}, we introduce notation used in this article and basic results of geometric measure theory and PDEs needed in the rest of the article. In section \ref{section.local.harnack.indecomposable}, we establish a criterion for indecomposability as well as a local Harnack inequality on balls and annuli around an indecomposable point. In section \ref{convergence.of.functions}, we introduce the notion of convergence of functions on converging multiplicity one varifolds. In section \ref{greens.function.stationary.varifolds}, we apply the Gr\"uter-Widman theory to construct the Green function on a stationary integral varifold. We also compute explicitly the Green function on a stationary cone.   In section \ref{convergence.of.green.function}, we show that the Green functions converge under multiplicity one convergence of varifolds with codimension 2 singular sets. In section \ref{behavior.green.function.singular.point}, we prove the asymptotic behavior of the Green function near an indecomposable point as well as global upper and lower bounds. In section \ref{remark.examples}, we discuss results of the Green function on stationary varifolds with decomposable points.
\section{Preliminaries}\label{preliminaries}

\subsection{Notation and basic concepts}
\label{notations}

We introduce the following notations:



$B_r(x)$
: the set $\{y\in \R^n:|x-y|<r\}$ (write $B_r$ if $x=0$);

$A_{t,r}(x)$: the set $B_{(1+t)r}(x)\setminus \overline{B_{r/(1+t)}(x)}$ for $x\in \R^n$, $t\in (0,1]$, $r>0$ (write $A_{t,r}$ if $x=0$);








$(E)_{\delta}$: the set $\{x\in \R^{n}: d(x, E)\le \delta\}$ for $E\subset \R^n$, $\delta>0$;

$\eta_{x,\lambda}$: the map $\R^n\to \R^n$ given by $\eta_{x,\lambda}(y)=\lambda(y-x)$ for $x\in \R^n$, $\lambda>0$ (write $\eta_\lambda$ if $x=0$);


$W^\perp$: the orthogonal complement of a linear subspace $W\subset \R^n$;






$\nabla f(x)$: the Euclidean gradient of $f$ at $x$;



$\omega_n$: the number $\HH^n(B_1)$;




$L^p_*(X,\mu)$: the space of weak $L^p$ functions on the measure space $(X,\mu)$;







$\SSS^n$: the standard $n$-sphere $\{x\in \R^{n+1}:|x|=1\}$;



$\Delta_\Sigma f$: the Laplace-Beltrami operator of $\Sigma$ acting on $f$;




The following are fundamental concepts in geometric measure theory. We refer to \cite{Fed,morgan2016gmt,Sim} for their precise definitions.

Let $n\ge 1,k\ge 0$, $U\subset\R^{n+k}$ be an open set. We introduce:

$\IV_n(U)$: the set of integral $n$-varifolds in $U$;


$T_xV$: the approximate tangent plane of $V\in \IV_n(U)$ at $x$;


$\nabla^V f(x)$: $\proj_{T_xV}\nabla f(x)$ for $f\in C^1(U)$;

$\nabla^\perp f(x)$: $\proj_{T_xV^\perp}\nabla f(x)$ for $f\in C^1(U)$;

$\|V\|$: the weight measure associated to $V\in \IV_n(U)$;

$V\llcorner W$, $\mu\llcorner W$: the restriction of $V\in \IV_n(U)$ or a Radon measure $\mu$ to a set $W\subset U$;

$|M|$: the varifold associated to a smooth submanifold $M^n\subset \R^{n+k}$;

$\phi_\#V$: the push forward of $V\in \IV_n(U)$ via a proper smooth map $\phi:U\to W$ where $W\subset \R^m$ is open;

$\spt V$, $\spt \chi$: the support of $V\in \IV_n(U)$ or $\chi\in C_c^\infty(U)$;

$\sing V$: the singular set of $V\in \IV_n(U)$;

$\reg V$: the regular set of $V\in \IV_n(U)$;

$\reg_\delta V$: the set $\{x\in \reg V:d(x,\sing V)> \delta\}=\spt V\setminus(\sing V)_\delta$ for $\delta>0$;

$\mathbf{F}_W(V_1,V_2)$: the varifold distance between two varifolds $V_1,V_2\in \IV_n(U)$ on $W\subset U$ see \cite[page 66]{pitts}.
\subsection{Stationary integral varifolds}
\label{Stationary.integral.varifolds}
In this section, we recall concepts and results related to integral varifolds. For more details, see \cite{allard1972first,Fed,Sim}. We assume $n\ge2$, $k\ge 1$ and $U\subset\R^{n+k}$ is an open set.

	Let $V\in \IV_n(U)$. For any $X\in C_c^1(U; \R^{n+k})$ the first variation of $V$ is \[\delta V(X)=\int_U \Div_{V}Xd\|V\|\]
	where $\Div_{V}X=\sum_{i=1}^{n} \nabla_{e_i}X\cdot e_i$ and $\{e_i\}_{i=1}^n$ is an orthonormal basis of $T_xV$. 
    We say that $V$ is stationary if  $\delta V=0$.

If $V=|M|$ is a smooth submanifold, then the first variation of $V$ involves the mean curvature in the interior and the boundary measure and the conormal. Now we recall the distributional boundary of a set $E$ in a varifold $V$ introduced in \cite{Me2weaklydifferentiable} as well as a coarea formula for rectifiable varifolds. The distributional boundary is to characterize the boundary of $E$ in $V$.
\begin{definition}[{\cite[Definition~5.1]{Me2weaklydifferentiable}}]\label{distributional.boundary}
    Let $V\in \IV_n(U)$ be a varifold such that $\|\delta V\|$ is a Radon measure and $E\subset \R^{n+k}$ be $\|V\|+\|\delta V\|$-measurable. We define the \textbf{distributional boundary} of $E$ in $V$ as
    \[V\partial E=(\delta V)\llcorner E-\delta (V\llcorner E).\]
\end{definition}
We recall the following coarea formula for rectifiable varifolds expressed in the language of distributional boundary. 
Note that \cite[Corollary~12.2]{Me2weaklydifferentiable} require $V$ to have $L^p$-generalized mean curvature (cf. \cite[Definition~16.5]{Sim}) to make sure weakly differentiable functions are $\|V\|$-approximately differentiable $\|V\|$-a.e. We do not need this assumption since we only consider Lipschitz functions and they are $\|V\|$-approximately differentiable $\|V\|$-a.e. by \cite[Lemma 4.5]{menne2012decay}.
\begin{theorem}[{\cite[Theorem 12.1, Corollary~12.2]{Me2weaklydifferentiable}}]\label{coarea}Let $V\in \IV_n(U)$, $f:U\to \R$ be a Lipschitz function and $E(y)=\{x\in U:f(x)>y\}$ for $y\in \R$. Then
    \[\int |\nabla^V f|d\|V\|=\int_0^\infty \|V\partial E(y)\|(U)dy.\]
    In particular, if $0\in U$, then for $\LL^1$-a.e. $0\le t<d(0,\partial U)$
    \[\frac{d}{dt} \|V\|(B_t)\ge \|V\partial B_t\|(U).\]
\end{theorem}

We also need the isoperimetric inequality from \cite[Theorem 7.1]{allard1972first}. The following version is not stated exactly as in \cite[Theorem 7.1]{allard1972first}, but its proof is contained in the proof of \cite[Theorem 7.1]{allard1972first} as long as one replace $V_\vp$ by $V\llcorner E$ in the proof. See also \cite[Theorem 2.2]{menne09isoperi} for a version where $V$ has $L^p$-generalized mean curvature. Note that the last inequality follows from Definition~\ref{distributional.boundary}.
\begin{theorem}[{\cite[Theorem 7.1]{allard1972first}}]\label{isoperimetric}
Let $V\in \IV_n(U)$ and $E\subset \R^{n+k}$ be $\|V\|+\|\delta V\|$-measurable. Then we have
\[\|V\|(E)^{(n-1)/n}\le \gamma(n,k)\|\delta (V\llcorner E)\|(U)\le \gamma(n,k)(\|\delta V\|(E)+\|V\partial E\|(U))\]
where $\gamma=\gamma(n,k)$ is some constant depending only on $n,k$.
\end{theorem}
The isoperimetric inequality and the coarea formula lead to the following lower bound on the mass if $V$ has controlled first variation cf. \cite[Lemma 1]{bombieri1972harnack}, \cite[8.3]{allard1972first}, \cite[2.5]{menne09isoperi}.
\begin{lemma}\label{volbd}
    Let $0<a<b<r$ and $V\in \IV(B_r)$ be such that $\|V\|(B_t)^{(n-1)/n}> 2\gamma\|\delta V\|(B_t)$ for all $t\in(a,b)$. Then $\|V\|(B_b)\ge c(b-a)^n$ where $c=c(n,k)=(2\gamma n)^{-n}$ and $\gamma$ is as in Lemma~\ref{isoperimetric}.
\end{lemma}
\begin{proof}
    Applying the isoperimetric inequalty Theorem~\ref{isoperimetric} and the coarea formula Theorem~\ref{coarea}, we have
    \[\|V\|(B_t)^{(n-1)/n}\le \gamma(\|\delta V\|(B_t)+\frac{d}{dt}\|V\|(B_t)).\]
    Using the assumption, we get
	\[\frac{d}{dt}\|V\|(B_t)\ge \frac{1}{2\gamma} \|V\|(B_t)^{(n-1)/n}.\]
    Integrating the inequality we get the desired result.
\end{proof}


By the monotonicity formula \cite[17.5]{Sim} for a stationary varifold, $r\mapsto \frac{\|V\|(B_r(x))}{\omega_nr^n}$ is non-decreasing. Thus for any $x\in \spt V$ the limit
\[\Theta(V,x)=\lim_{r\to0}\frac{\|V\|(B_r(x))}{\omega_nr^n}\]
exists. A varifold $V\in \IV_n(\R^{n+k})$ has Euclidean volume growth if there exists $\Lambda\ge 1$ such that 
\begin{equation}\label{euclidean.volume.growth}
	\|V\|(B_R(x))\le \Lambda \omega_n R^n\text{ for any }x\in \spt V\text{ and }R>0.
\end{equation}Note that when we say that a sequence $\{V_i\}\subset \IV_n(\R^{n+k})$ has Euclidean volume growth \eqref{euclidean.volume.growth}, we mean that \eqref{euclidean.volume.growth} holds for all $V_i$ with the same $\Lambda$.
If $V$ has Euclidean volume growth, then the following limit
\[\Theta(V,\infty)=\lim_{r\to\infty}\frac{\|V\|(B_r)}{\omega_nr^n}\] also exists
and we have $\Theta (V,x)\le \Theta(V,\infty)$ for any $x\in \spt V$.

We say that a stationary varifold $V\in \IV_n(U)$ has multiplicity one if
\begin{equation}\label{mult1}
	\Theta(V, x)=1\quad \forall x\in \reg V.
\end{equation}

	Let $\C\in \IV_n(\R^{n+k})$. We say that $\C$ is a cone if $\eta_{\lambda\#}\C=\C$ for any $\lambda>0$.
	
	Let $V\in \IV_n(U)$ be stationary. We say that $\C\in \IV_n(\R^{n+k})$ is a tangent cone of $V$ at $x$ (resp. at infinity) if there is a sequence $\lambda_j\to0$ (resp. $\lambda_j\to \infty$) such that $\eta_{x,\lambda_{j}\#}V\rightharpoonup \C$ in $\IV_n(\R^{n+k})$. We write $\tan(V,x)$ (resp. $\tan(V,\infty)$) for the set of all tangent cones of $V$ at $x$ (resp. at infinity). It is a consequence of the monotonicity formula \cite[17.5]{Sim} that tangent cones are cones.

\begin{remark}\label{known.result.sing.v}
	Very little is known about the singular set of a general stationary integral varifold. By \cite[8.1(1)]{allard1972first}, we know that $\sing V$ is meager i.e. its complement $\reg V$ is an open dense subset of $\spt V$. 
	
	If $V$ is the associated varifold of an area-minimizing current in any codimension, then $\HH^{n-2+\alpha} (\sing V)=0$ for any $\alpha>0$ (see e.g. \cite[Theorem 3.4]{DeLellis2016}). If $V$ is the associated varifold of an area-minimizing boundary, then $V$  satisfies \eqref{euclidean.volume.growth} with $\Lambda=\Lambda(n)$ by \cite[Proof of 37.2 (1)]{Sim}, $V$ satisfies (\ref{mult1}) by \cite[Remark 37.1 (**)]{Sim} and $\HH^{n-7+\alpha} (\sing V)=0$ for any $\alpha>0$ by \cite{Fed70}.
\end{remark}
We need the following capacity argument to get a sequence of cut-off function which will be useful in Lemma~\ref{getting.sobolev.functions}.
\begin{lemma}\label{cutoff.function}
	Let $1\le p<n$, $W\ssubset U$, $S\subset U$ and $V\in \IV_n(U)$ be such that $\|V\|(B_r(x))\le C(n,W,V)r^n$, for any $x\in W$ and $r\le \frac{1}{2}d(W,\partial U)$. If $1<p<n$, we assume $\HH^{n-p}(S)<\infty$. If $p=1$, we assume $\HH^{n-1}(S)=0$. Then for $0<\delta<\frac{1}{2}d(W,\partial U)$, there exists $\{\chi_i\}_{i=1}^\infty\subset C_c^\infty(\R^{n+k})$ such that \begin{enumerate}
		\item [(i)] $0\le \chi_i\le 1$,
		\item[(ii)]$S\cap \bar W\subset \{\chi_i=1\}^\circ$,
		\item[(iii)]$\chi_i=0$ on $\R^{n+k}\setminus (S\cap \bar W)_\delta$,
		\item[(iv)] $\chi_i\to 0$ pointwise on $\R^{n+k}\setminus(S\cap \bar W)$,
		\item[(v)]$\int_U |\nabla^V \chi_i|^pd\|V\|\to0$ as $i\to\infty$.
	\end{enumerate}
\end{lemma}
\begin{proof}
	The proof is a slight modification of the standard proof showing the $p$-capacity of $S$ is 0 cf. \cite[Theorem 4.16]{gariepy2015measure}. If $p=1$, then $\chi_i$ is constructed as $f$ in the proof of \cite[Theorem 4.16]{gariepy2015measure}. If $1<p<n$, $\chi_i$ is constructed exactly as the $g_i$ in the proof of \cite[Theorem 4.16]{gariepy2015measure}. The differences are as follows. First, when $p=1$, since $\HH^{n-1}(S)=0$, for any $i\in\N$ we take a cover $\{B^0(x_j,r_j)\}_{j=1}^m$ such that $\sum r_j^{n-1}<i^{-1}$. Moreover, we note that the integrand in (v) is $|\nabla^V\chi_i|^p$ instead of $|\nabla \chi_i|^p$ and that the integration is with respect to $\|V\|$. The first issue can be addressed by the fact that $|\nabla ^Vf|\le |\nabla f|$ for any Lipschitz function $f$. The second issue is addressed by the assumption $\|V\|(B_r(x))\le Cr^n$. Hence the estimates of $\int_U |\nabla^V\chi_i|^pd\|V\|$ in the proof of \cite[Theorem 4.16]{gariepy2015measure} work in our setting.
\end{proof}

\subsection{Preliminaries on Sobolev spaces and PDEs on stationary integral varifolds}
In this section, we introduce the theory of Sobolev spaces and elliptic PDEs on stationary integral varifolds. We note that many of the results in subsections \ref{section.sobolev.space}, \ref{section.poisson.equation} hold for more general varifolds e.g. varifolds with bounded generalized mean curvature, cf. \cite{Me1sobolev,Me2weaklydifferentiable,MS}. We stick to the stationary case for simplicity.
\subsubsection{Sobolev spaces on stationary integral varifolds}\label{section.sobolev.space}

\begin{definition}[\cite{Me1sobolev}]\label{Sobolev.space}Let $1\le p\le\infty$ and $V\in \IV_n(U)$ be stationary and $U_1\subset U$ be an open subset.
	We say that a function $f\in L^p_\loc(U_1,\|V\|)$ is in $W^{1,p}_\loc(U_1,V)$ if there exists $F\in L^p_\loc (U_1,\|V\|; \R^{n+k})$ so that for any $\epsilon>0$ and $K\subset U_1$ compact, there is $g\in \lip(U_1)$ so that \[\|f-g\|_{L^p(K, \|V\|)}+\|F-\nabla^Vg\|_{L^p(K, \|V\|)}<\epsilon.\]
    We write $\nabla^Vf$ for $F$.
	The Sobolev norm is defined as \[\|f\|_{W^{1,p}(U_1,V)}=(\|f\|_{L^p(U_1,\|V\|)}^p+\|\nabla^Vf\|_{L^p(U_1,\|V\|)}^p)^{1/p}.\]
	
	The Sobolev space $W^{1,p}(U_1,V)$ consists of $f\in W^{1,p}_\loc(U_1,V)$ with $\|f\|_{W^{1,p}(U_1,V)}<\infty$.
	
	The Sobolev space with zero boundary value $W^{1,p}_0(U_1,V)$ is the closure of $\lip(U_1)$ with compact support in $U_1$ under the norm $\|\cdot\|_{W^{1,p}(U_1,V)}$.

    If $U_1=U$, we write $W^{1,p}_\loc(V)$, $W^{1,p}(V)$, $W^{1,p}_0(V)$ for $W^{1,p}_\loc(U,V)$, $W^{1,p}(U,V)$, $W^{1,p}_0(U,V)$  respectively.
\end{definition}

\begin{remark}\label{remark.approx.ambient.functions}
	If $g\in \lip (U)$, by \cite[Example 8.7]{Me2weaklydifferentiable}, $g$ is $\|V\|$-approximately differentiable $\|V\|$-a.e. in $U$. At such a point, $\nabla^V g(x)=\proj_{T_xV}(\nabla g(x))$ where $\nabla g$ is the approximate differential of $g$ at $x$.
\end{remark}

The Sobolev functions on varifolds share many properties with those on Eulidean spaces.  We summarize those which are useful in later applications.
For Sobolev functions, the integration by parts formula holds if we compose the function with a smooth function with compact support.

\begin{proposition}[{\cite[Remark 5.2]{Me1sobolev}}]\label{integration.by.parts}
	Let $V\in \IV_n(U)$ be stationary, $X\in C_c^1(U;\R^{n+k})$, $\gamma\in C_c^\infty(\R)$ and $f\in W^{1,1}_\loc(V)$. Then 
	\[\int \nabla (\gamma\circ f)\cdot Xd\|V\|=-\int (\gamma\circ f)\Div_{V}X d\|V\|.\]
\end{proposition}
\begin{remark}
	Note that if $|f|\le C$ on $\spt V$, then we can take $\gamma\in C_c^\infty(\R)$ so that $\gamma(t)=t$ on $[-C,C]$. Then the integration by parts formula is true for $f$. This will always be the case whenever we apply Proposition~\ref{integration.by.parts}.
\end{remark}

In the following two lemmas, we show that a Sobolev function can be extended accross a set $S$ if we control the behavior of the function near $S$ and the size of $S$. A typical choice is $S=\sing V$. We recall that we write $(S)_\delta=\{x\in \R^{n+k}:d(x,S)\le\delta\}$.

\begin{lemma}\label{sobolev.function.vanishing.on.singular.set}
	Let $1\le p<\infty$ and $V\in \IV_n(U)$ be stationary. 
    Let $S\subset U$ be a closed set and $f:\spt V\to \R$ be such that $f\in W^{1,p}_\loc(U\setminus S,V)$. If $f=0$ on $\spt V\cap (S)_{2\delta}$ for some $\delta>0$, then $f\in W^{1,p}_\loc(V)$.
\end{lemma}
\begin{proof}
	Let $W\ssubset U$ be fixed. 
    There exist $f_i\in \lip(U)$ such that $f_i\to f$ in $W^{1,p}(W\setminus (S)_{\delta},V)$. We take $0\le \zeta_\delta\in C_c^{\infty}(\R^{n+k})$ such that $\zeta_\delta=1$ on $\R^{n+k}\setminus (S\cap \bar W)_{2\delta}$ and $\zeta_\delta =0$ on $(S\cap \bar W)_{\delta}$. Then using $|\zeta_\delta|\le 1$, $\spt |\nabla\zeta_\delta|\subset (S\cap \bar W)_{2\delta}$ and $f=0$ on $(S\cap \bar W)_{2\delta}$ we have
	\[\int_W |\zeta_\delta f-\zeta_\delta f_i|^pd\|V\|+\int_W |\nabla(\zeta_\delta f-\zeta_\delta f_i)|^pd\|V\|\le \|f-f_i\|_{W^{1,p}(W\setminus (S)_{\delta},V)}^p\to0.\] 
	Thus $\zeta_{\delta}f\in W^{1,p}(W,V)$. Since $f=0$ on $(S\cap \bar W)_{2\delta}$, we have $f=\zeta_\delta f\in W^{1,p}(W,V)$.
\end{proof}

\begin{lemma}\label{getting.sobolev.functions}
	Let $1< p<n$, $V\in \IV_n(U)$ be stationary and $S\subset U$ be closed and such that $\HH^{n-p}(S)<\infty$. 
    If $f\in W^{1,p}_\loc( U\setminus S,V)\cap L^\infty((S)_{\delta}\cap U,\|V\|)$ for some $\delta>0$ and $\nabla ^V f\in L^p_\loc(U,\|V\|)$, then $f\in W^{1,p}_\loc(V)$.
\end{lemma}

\begin{proof}
	Let $W\ssubset U$ be fixed. Since $\HH^{n-p}(S)<\infty$, let $\chi_i$ be the cut-off function for $S\cap \bar W$ constructed in Lemma~\ref{cutoff.function}. We have $\spt |\nabla\chi_i|\subset (S\cap \bar W)_{\delta}$, $\chi_i\to 0$ pointwise on $\R^{n+k}\setminus(S\cap \bar W)$ and $\int_U |\nabla^V \chi_i|^pd\|V\|\to0$. By Lemma~\ref{sobolev.function.vanishing.on.singular.set}, $(1-\chi_i) f\in W^{1,p}( W,V)$ for any $i\in \N$. We claim that $(1-\chi_i)f\to f$ in $W^{1,p}(W,V)$. Indeed, we have
	\[\begin{aligned}
		&\int_W |\chi_i f|^pd\|V\|+\int_W |\nabla^V(\chi_i f)|^p d\|V\|\\
		\le& \int_W \chi_i^p f^pd\|V\|+C\int_W \chi_i^p|\nabla^V f|^pd\|V\|+C\left(\sup_{(S)_{\delta}\cap W} |f|^p\right)\int_{(S)_{\delta}\cap W} |\nabla^V\chi_i|^p d\|V\|.
	\end{aligned}\]
	Since $\HH^{n-p}(S)<\infty$, we have $f\in L^p(W,\|V\|)$. By $L^\infty((S)_{\delta}\cap W,\|V\|)$ and $\nabla^V f\in L^p(W,\|V\|)$, the above terms tend to $0$ as $i\to\infty$ by Lebesgue dominated convergence.
\end{proof}

For stationary integral varifolds, we have the Sobolev inequality similar to that on $\R^n$.

\begin{theorem}[\cite{allard1972first}, \cite {MS}, {\cite[Theorem 7.2]{Me1sobolev}}]\label{Michael.Simon.Sobolev.inequality}
	Let $1\le p<n$ and $V\in \IV_n(U)$ be stationary. For any $f\in W^{1,p}_0(V)$, we have
	\[\|f\|_{L^{np/(n-p)}(U,\|V\|)}\le C(n)\|\nabla ^Vf\|_{L^p(U,\|V\|)}.\]
	For any $W\ssubset U$ and $f\in W^{1,p}(V)$, we have
	\[\|f\|_{L^{np/(n-p)}(W,\|V\|)}\le C(n,W)\|f\|_{W^{1,p}(V)}.\]
\end{theorem}


We also have the Rellich compactness theorem for Sobolev spaces on stationary varifolds.
\begin{theorem}[{\cite[Theorem 7.11, 7.21]{Me1sobolev}}] \label{rellich.on.one.varifold}
	Let $1\le p<n$ and $V\in \IV_n(U)$ be stationary with $n\ge 2$. Then bounded subsets in $W^{1,p}_\loc(V)$ (resp. $W^{1,p}_0(V)$) have compact closure in $L^q_\loc(U,\|V\|)$ (resp. $L^q(U,\|V\|)$) where $1\le q<\frac{np}{n-p}$.
\end{theorem}

\subsubsection{Poisson equation on stationary integral varifolds}\label{section.poisson.equation}
\begin{definition}
	Let $V\in \IV_n(U)$ be stationary. We say that a function $u\in W^{1,2}_\loc(V)$ solves the equation $-\Delta_V u=f$ in $U$ if 
	\begin{equation}\label{weaksln}
		\int \nabla^V u\cdot \nabla^V \phi d\|V\|=\int f\phi d\|V\|\quad \forall\phi\in C_c^\infty(U).
	\end{equation} 
	We also say that $-\Delta_V u\le f$ in $U$ if we have $\le$ sign in (\ref{weaksln}) for all  $0\le \phi\in C_c^\infty(U)$.
\end{definition}


\begin{lemma}[Caccioppoli inequality]\label{cac}
	Let $V\in \IV_n(U)$ be stationary. Let $u\in W^{1,2}_\loc(V)$ solve $-\Delta_V u=f$ in $U$ for $f\in L^2(U,\|V\|)$. 
	Then for any $R>0$ such that $ B_{2R}(x)\ssubset U$, we have
	\[\int_{B_R(x)} |\nabla^V u|^2d\|V\|\le C\left(\frac{1}{R^2}\int_{B_{2R}(x)}u^2d\|V\|+R^2\int_{B_{2R}(x)} f^2d\|V\|\right)\]
	where $C$ is an absolute constant.
\end{lemma}
\begin{proof}
	The conclusion follows from plugging $\phi=\psi^2u$ in (\ref{weaksln}) with $\psi\in C_c^\infty(B_{2R}(x))$, $0\le \psi\le 1$, $\psi=1$ on $B_R(x)$ and $|\nabla\psi|\le C/R$ and the applying Cauchy-Schwarz inequality.
\end{proof}
We have the following mean value inequality for subharmonic functions.

\begin{theorem}[Mean value inequality {\cite[7.5 (6)]{allard1972first}, \cite[Theorem 3.4]{MS}}]\label{mean.value.inequality}
	Let $V\in \IV_n(U)$ be stationary. For any $B_R(x)\ssubset U$, $0\le u\in W^{1,2}(B_R(x),V)$, $\Delta_Vu\ge 0$ in $B_R(x)$ we have
	\[\underset{\spt V\cap B_{R/2}(x)}{\esssup} u\le C(n) R^{-n}\int_{B_R(x)}ud\|V\|.\]
\end{theorem}
Next, we need an \textit{a priori} estimate for subsolutions to Poisson equation. This also gives a weak maximum principle for subharmonic functions on stationary varifolds.

\begin{theorem}[\textit{A priori} estimate]\label{apri}Let $V\in \IV_n(U)$ be stationary and $W\ssubset U$. Let $u\in W^{1,2}(W,V)$ satisfy $-\Delta_V u\le f$ for $f\in L^{q}(W,\|V\|)$ with $q>n/2$. We define $\sup\limits_{V\partial W}u$ as the minimal $l$ such that $(u-l)_+\in W^{1,2}_0(W,V)$. Then there exists $C=C(n,k,q,\|V\|(W))$ such that 
	\[\underset{\spt V\cap W}{\esssup}\ u\le \sup_{V\partial W}u^++C\|f\|_{L^{q}(W,\|V\|)}.\]
\end{theorem}
\begin{proof}
	This theorem follows from the same argument as \cite[Theorem 8.15 and 8.16]{GT}. In applying the argument, we are able to take the test functions due to \cite[Remark 5.6]{Me1sobolev}. We also need to replace the Euclidean Sobolev inequality by Theorem~\ref{Michael.Simon.Sobolev.inequality}.
\end{proof}

\section{Local Harnack inequality around an indecomposable point}
\label{section.local.harnack.indecomposable}
\subsection{A criterion for indecomposability}
We start by proving the following characterization of indecomposability of a varifold by the multiplicity and connectivity of its regular set. Recall Definition~\ref{definition.indecomposable} that $V$ is indecomposable if $V=V_1+V_2$ and $\delta V_1=\delta V_2=0$ implies either $V_1=0$ or $V_2=0$.
 
\begin{lemma}    \label{reg.connected.mult1.equiv.indecomposable}
    Let $V\in \IV_n(U)$ be stationary. 
    
    (i) If $\mathcal{H}^n(\sing V) = 0$, $V$ has multiplicity one and $\reg V$ is connected, then $V$ is indecomposable.  

    (ii) Conversely, if $\mathcal{H}^{n-1}(\sing V) = 0$ and $V$ is indecomposable, then $V$ has multiplicity one and $\reg V$ is connected.

\end{lemma}
\begin{proof}(i) Let $V=V_1+V_2$ with $\delta V_1=\delta V_2=0$ and $x\in \reg V$. We show that $x\in \reg V_1\cup \reg V_2$. We have $\Theta (V_i,x)=0$ or $\Theta (V_i,x)\ge 1$ by the monotonicity formula for $V_i$, $i=1,2$.  Since $V=V_1+V_2$, $1=\Theta(V,x)=\Theta (V_1,x)+\Theta (V_2,x)$. Thus either $\Theta(V_1,x)=1$ or $\Theta(V_2,x)=1$. By Allard's regularity theorem \cite[Theorem 24.2]{Sim}, $x\in \reg V_1$ or $x\in \reg V_2$ and $\reg V_1\cap \reg V_2 \cap \reg V=\emptyset$. Therefore, $\reg V=(\reg V_1\cap \reg V)\sqcup (\reg V_2\cap \reg V)$. Since $\reg V$ is connected we have either $\reg V_1\cap \reg V=\emptyset$ or $\reg V_2\cap \reg V=\emptyset$. By $\HH^{n}(\sing V)=0$, we have $V_1=0$ or $V_2=0$.
    
    (ii) Let $M$ be a connected component of $\reg V$. Define $V_1=|M|$. Since the mean curvature $H=0$ on $M$, we have
    \begin{align*}
        \delta V_1(X)=&\int_M \Div_M (X)d\|V\|\\
        =&-\int_M \langle H,X^\perp\rangle d\|V_1\| +\int _M \Div _M(X^\top)d\|V_1\|\\
        =&\int _M \Div _M(X^\top)d\|V_1\|.
    \end{align*}
    We restrict ourself to a bounded open set $U_1$ containing $\spt X$. Since $\HH^{n-1}(\sing V)=0$ and $V$ is stationary, there exists $\chi_i\in C_c^\infty(U_1)$ satisfying Lemma~\ref{cutoff.function}. By Sard's theorem, coarea formula and the Markov inequality, we take $M_i=\{x\in M\cap U_1:\chi_i(x)\le c_i\}$ where $c_i$ is a regular value of $\chi_i$ such that $\HH^{n-1}(\partial M_i)\le 2\int |\nabla^V\chi_i|d\|V\|\to0$. We can arrange $c_i$ and $\chi_i$ such that $M_i\subset M_{i+1}\subset M\cap U_1$ and $\cup_{i=1}^\infty M_i=M\cap U_1$. We apply the divergence theorem on $M_i$ and take $i\to\infty$ to get $\delta V_1=0$. Thus $V_1=|M|\in \IV_n(U)$ is a stationary. Since $M$ is a connected component of $\reg V$, there is $V_2\in \IV_n(U)$ such that $V=V_1+V_2$. Since $\delta V=\delta V_1=0$, we have $\delta V_2=0$. By the indecomposability of $V$, we have $V_2=0$ and $V=V_1=|M|$. Thus $\reg V=M$ is connected and $V=|\reg V|$ has multiplicity one.
\end{proof}
We now present an example showing that the assumption $\HH^{n-1}(\sing V)=0$ in (ii) is sharp.
\begin{lemma}\label{triple.junction.indecomposable}
    Let $H_1,H_2,H_3$ be three half $n$-planes in $\R^{n+k}$ with $\partial H_1=\partial H_2=\partial H_3=\R^{n-1}$ and $H_i$ and $H_j$ have angle $120^\circ$ for $i\neq j$. Let $V=H_1+H_2+H_3$. It is well-known that $V$ is stationary. Then $V$ is indecomposable.
\end{lemma}

\begin{proof}
    Suppose $V=V_1+V_2$ where $V_i\in \IV_n(\R^{n+k})$ and $\delta V_i=0$ for $i=1,2$. As in the proof of Lemma \ref{reg.connected.mult1.equiv.indecomposable}, we have $\reg V=(\reg V_1\cap \reg V)\sqcup(\reg V_2\cap \reg V)$. Since $\reg V=H_1\cup H_2\cup H_3$ has 3 components, for any $i=1,2,3$ there is $j=1,2$ such that $H_i\subset \reg V_j$. Thus $V_1,V_2$ are unions of half planes. Since $H_i$ is not stationary, the only way for both $V_1$ and $V_2$ to be stationary is if $V_1=0$ or $V_2=0$.
\end{proof}
\begin{remark}\label{area.min.bdy.indecomposable}
    If $V$ is the associated varifold of an area-minimizing boundary, then $V$ is indecomposable by \cite[Theorem 1]{bombieri1972harnack}. In particular,  every point in $\spt V$ as well as infinity is indecomposable.
\end{remark}
\begin{remark}\label{different.notion.indecomposable}
    There are two other definitions of indecomposability of a varifold in the literature, namely \cite[Definition~2.15]{Modino} and \cite[Definition~6.2]{Me2weaklydifferentiable}. When restricted to stationary varifolds, they have the following relation
    \[\text{\cite[Definition~2.15]{Modino}} \Longrightarrow \text{Definition~\ref{definition.indecomposable}}\Longrightarrow \text{\cite[Definition~6.2]{Me2weaklydifferentiable}}.\]
    It is known that Definition~\ref{definition.indecomposable} is more restrictive than \cite[Definition~6.2]{Me2weaklydifferentiable} since for any $V$ indecomposable, Definition~\ref{definition.indecomposable} excludes $kV$ for $k\in \Z$, $|k|\ge 2$ to be indecomposable while it is still indecomposable in the sense of \cite[Definition~6.2]{Me2weaklydifferentiable}.

    It is unknown whether \cite[Definition~2.15]{Modino} is equivalent to Definition~\ref{definition.indecomposable} as \cite[Definition~2.15]{Modino} allow $V_1,V_2$ to have $L^p$-generalized mean curvature.
\end{remark}

    



    

\subsection{Local Harnack inequality on balls}
The following theorem shows that near an indecomposable point of $V$, any stationary varifold that is close to $V$ supports a relative isoperimetric inequality in a neighborhood of the point. In this way, the indecomposability of a point propagates to a neighborhood of the point.
\begin{theorem}[Relative isoperimetric inequality]\label{relative.isoperimetric.varifold}
	    Let $V\in \IV_n(U)$ be stationary and $x\in\indec V$. Then there exists $\beta_0=\beta_0(n,V,x)\in(0,1)$, $r_0=r_0(n,V,x)>0$, $\epsilon_0=\epsilon_0(n,V,x)\in(0,1)$ with the following significance.  For any $r<r_0$, $V'\in \IV_n(B_r(x))$ stationary with $\mathbf{F}_{B_r(x)}(V,V')<\epsilon_0 r^n$ and $\|V'\|$-measurable set $E\subset \R^{n+k}$, we have
	\[\|V'\partial E\|(B_r(x))=\|V'\partial E^c\|(B_r(x))\ge \frac{1}{2\gamma}\min\{\|V'\|( E\cap B_{\beta_0 r}(x)),\|V'\|(E^c\cap B_{\beta_0 r}(x))\}^{\frac{n-1}{n}}\]
	where $\gamma$ is the constant from Lemma~\ref{isoperimetric}.
\end{theorem}

\begin{proof}Without loss of generality, we assume $x=0$.
	Suppose the conclusion is not true. Then for any $j\in\N$, $\beta_0=r_0=\epsilon_0=j^{-1}$, there exists $r_{j}<j^{-1}$, $V_j\in \IV(B_{r_j})$ stationary  with $\mathbf{F}_{B_{r_j}}(V_j,V)<j^{-1} r_j^n$ and $E_j\subset \R^{n+k}$ such that \[\|V_j\partial E_j\|(B_{r_j})=\|V_j\partial E_j^c\|(B_{r_j})< \frac{1}{2\gamma}\min\{\|V_j\|( E_j\cap B_{j^{-1} r_j}),\|V_j\|(E_j^c\cap B_{j^{-1} r_j})\}^{\frac{n-1}{n}}.\]
    We put $V_j^{(1)}=V_j\llcorner E_j\in \IV_n(B_{r_{j}})$ and $V_j^{(2)}=V_j\llcorner E_j^c\in \IV_n(B_{r_{j}})$ so that $V_j=V_j^{(1)}+V_j^{(2)}$. Since $V_j$ is stationary, we have $\delta V_j^{(1)}=-V_j\partial E_j$ and $\delta V_j^{(2)}=-V_j\partial E_j^c$. Thus
    \[\|\delta V_j^{(1)}\|(B_{r_j})=\|\delta V_j^{(2)}\|(B_{r_j})<\frac{1}{2\gamma}\min\{\|V_j^{(1)}\|( B_{j^{-1} r_j}),\|V_j^{(2)}\|(B_{j^{-1} r_j})\}^{\frac{n-1}{n}}.\]
     
	We consider the rescaling $\tilde V_j=\eta_{r_j^{-1}\#}(V_j)$ and $\tilde V_j^{(i)}=\eta_{r_j^{-1}\#}(V_j^{(i)})$, $i=1,2$. Then we have $\tilde V_j=\tilde V_j^{(1)}+\tilde V_j^{(2)}$ and
	\begin{equation}\label{delta V<min V^(n-1)/n.rescaled}
	    \|\delta \tilde V_j^{(1)}\|(B_1)=\|\delta \tilde V_j^{(2)}\|(B_1)< \frac{1}{2\gamma}\min\{\|\tilde V_j^{(1)}\|(B_{j^{-1}}),\|\tilde V_j^{(2)}\|( B_{j^{-1}})\}^{\frac{n-1}{n}}.
	\end{equation}
	In particular for any $\frac{1}{j}\le t\le 1$
	\[
	    \|\delta \tilde V_j^{(i)}\|(B_t)< \frac{1}{2\gamma}\|\tilde V_j^{(i)}\|(B_t)^{\frac{n-1}{n}},\quad i=1,2.
	\]
	Thus by Lemma~\ref{volbd}, we have for $\frac{1}{j}\le t\le 1$
	\begin{equation}\label{mass.lower.bound}\|\tilde V_j^{(i)}\|(B_t)\ge c(t-\frac{1}{j})^{n},\quad i=1,2.\end{equation}
	Since $\tilde V_j^{(1)},\tilde V_j^{(2)}$ have locally bounded mass and first variation, by \cite[Theorem 42.7]{Sim} there is a subsequence still denoted by $\{j\}$ such that $\tilde V_j^{(i)}\to \C^{(i)}$ in $\IV_n(B_1)$ for $i=1,2$ and $\eta_{r_j^{-1}\#}V\to \C$ in $\IV_n(\R^{n+k})$ for some $\C\in \tan(V,0)$. Since $\mathbf{F}_{B_1}(\tilde V_j,\eta_{r_j^{-1}\#}V)\le r_j^{-n}\mathbf{F}_{B_{r_j}}(V_j,V)<j^{-1}$ and $\mathbf{F}_{B_1}(\eta_{r_j^{-1}\#}V,\C)\to0$, we have $\tilde V_j\to \C\llcorner B_1$ in $\IV(B_1)$. 
    
    Thus $\C\llcorner B_1=\C^{(1)}+\C^{(2)}$. For $i=1,2$, we have $\delta \C^{(i)}=0$ by \eqref{delta V<min V^(n-1)/n.rescaled} and $\|\C^{(i)}\|(B_r)\ge cr^n$ for $0\le r\le 1$ by \eqref{mass.lower.bound}. We consider a further rescaling $\C\llcorner B_j=\eta_{j\#}(\C \llcorner B_1)=\eta_{j\#}\C^{(1)}+\eta_{j\#}\C^{(2)}$. By \cite[Theorem 42.7]{Sim}, after taking a subsequence we have  $\eta_{j\#}\C^{(i)}\to \C_i$ in $\IV_n(\R^{n+k})$ as $j\to\infty$ and $\|\C_i\|(B_r)\ge cr^n$. Then $\C=\C_1+\C_2$ with $\delta \C_i=0$ and $\C_i\neq0$ for $i=1,2$.
	This contradicts the indecomposability of $0\in\spt V$.
\end{proof}
\begin{theorem}[Neumann-Sobolev inequality]\label{Neumann.sobolev.inequality}
    Let $V\in \IV_n(U)$ be a stationary and $x\in\indec V$. With the same $\beta_0,r_0,\epsilon_0,r,V'$ as in Theorem~\ref{relative.isoperimetric.varifold}, for any $f\in C^1(B_{r}(x))$, we have
    \[\min_{c\in \R}(\int _{B_{\beta_0 r}(x)}|f-c|^{\frac{n}{n-1}}d\|V'\|)^{\frac{n-1}{n}}\le 2\gamma \int_{B_{r}(x)}|\nabla^{V'} f|d\|V'\|.\]
\end{theorem}
\begin{proof}
    The proof follows from the same argument as \cite[Theorem 3]{bombieri1972harnack} with Theorem~\ref{relative.isoperimetric.varifold} replacing \cite[Theorem 2]{bombieri1972harnack}.
\end{proof}
\begin{theorem}[Local Harnack inequality]\label{LocalHarnack.balls}
	Let $V\in \IV_n(U)$ be stationary and $x\in \indec V$. Let $\beta_0,r_0,\epsilon_0,r,V'$ be as in Theorem~\ref{relative.isoperimetric.varifold}, and assume that $\spt V'\cap B_{\beta_0 r/2}(x)\neq\emptyset$. Then for any $0\le u\in W^{1,2}(B_{r}(x),V')$ satisfying $-\Delta_{V'}u=f$ in $B_{r}(x)$ for $f\in L^{q}(B_r(x),\|V'\|)$ with $q>n/2$, we have 
    \begin{equation}\label{harnack}
        \sup_{\spt V'\cap B_{\beta_0 r}(x)}u\le C(\inf_{\spt V'\cap B_{\beta_0 r}(x)} u+r^{2-\frac{n}{q}}\|f\|_{L^{q}(B_r(x),\|V'\|)})
    \end{equation}
    for some $C=C(n,k,q,\Lambda,\beta_0)>0$ where $\Lambda=\omega_n^{-1}r_0^{-n}\|V'\|(B_{r_0}(x))$.
	Moreover, $u\in C^\alpha(\spt V'\cap B_{\beta_0 r}(x))$ for some $\alpha=\alpha(n,k,q,\Lambda,\beta_0)\in (0,1)$ and  there exists $C=C(n,k,q,\Lambda,\beta_0)$ such that
	\[[u]_{\alpha,\spt V'\cap B_{\beta_0 r}(x)}\le Cr^{-\alpha}(\|u\|_{L^\infty(B_r(x),\|V'\|)}+r^{2-\frac{n}{q}}\|f\|_{L^{q}(B_r(x),\|V'\|)})\]
	where \[[u]_{\alpha,\spt V'\cap B_{\beta_0 r}(x)}:=\sup_{x\neq y\in \spt V'\cap B_{\beta_0 r}(x)}\frac{|u(x)-u(y)|}{|x-y|^\alpha}.\]
	
	If instead $0\le u\in W^{1,2} ( B_r(x),V')$ satisfies $\Delta_{V'} u\le0$ in $B_r(x)$, then for any $0<p<\frac{n}{n-2}$, there exists $C=C(n,k,q,\Lambda,\beta_0,p)>0$ such that 
	\[(\frac{1}{\|V'\|(B_{\beta_0 r}(x))}\int_{B_{\beta_0 r}(x)} u^pd\|V'\|)^{1/p}\le C\inf_{\spt V'\cap B_{\beta_0 r}(x)} u.\]
\end{theorem}
\begin{proof}
	In \cite{bombieri1972harnack}, Bombieri and Giusti proved the Harnack inequality for the $f=0$ case. One can modify the proof following \cite[Theorem 8.17 and Theorem 8.18]{GT} to show the Harnack inequality we have here. The necessary modifications are to use the Sobolev inequality from Theorem~\ref{Michael.Simon.Sobolev.inequality} and the Neumann-Sobolev inequality from Theorem~\ref{Neumann.sobolev.inequality} as well as the abstract John-Nirenberg inequality \cite[Theorem 4]{bombieri1972harnack}. The assumption $\spt V'\cap B_{\beta_0 r/2}(x)\neq\emptyset$ is to guarantee $\|V'\|(B_{\beta_0 r}(x))\ge \|V'\|(B_{\beta_0 r/2}(y))\ge \omega_n (\beta_0 r/2)^n$ where $y\in \spt V'\cap B_{\beta_0 r/2}(x)$. The $C^\alpha$ estimates follows from the same proof as in \cite[Theorem 8.24]{GT} once we have the Harnack inequality.
\end{proof}
\begin{remark}\label{area.min.harnack}
    If $V$ is the associated varifold of an area-minimizing boundary, then  Theorem~\ref{LocalHarnack.balls} is just \cite[Theorem 5]{bombieri1972harnack} and we have $r_0(n,V,x)=\infty$  and $\beta_0,\Lambda$ only depend on $n$. That $r_0(n,V,x)=\infty$ and that $\beta_0$ depends only on $n$ follows from the fact that any area-minimizing boundary is indecomposable. One can therefore extract a convergent subsequence from any contradicting sequence---without requiring it to arise from a blow-up of a fixed varifold---and obtain a contradiction to indecomposability. See \cite[Theorem~2]{bombieri1972harnack}. If $V$ is the associated varifold of a stable minimal hypersurface in $B_4$ away from multiplicity 2 varifolds, then $\spt V\cap B_1\subset \indec V$ by \cite[Proposition~3.1]{Wang} and Theorem~\ref{LocalHarnack.balls} is \cite[Corollary~3.5]{Wang}.
\end{remark}
\begin{corollary}\label{cone.r_0=infty}
    Let $\C\in \IV_n(\R^{n+k})$ be a stationary cone which is indecomposable at $0\in \spt \C$. Then Theorem~\ref{relative.isoperimetric.varifold}, \ref{Neumann.sobolev.inequality} and \ref{LocalHarnack.balls} hold with $r_0(n,\C,0)=\infty$.
\end{corollary}
\begin{proof}It suffices to show $r_0(n,\C,0)=\infty$ in Theorem~\ref{relative.isoperimetric.varifold} since the rest are consequences of Theorem~\ref{relative.isoperimetric.varifold}. By Theorem~\ref{relative.isoperimetric.varifold}, there are $\beta_0,r_0,\epsilon_0$ such that we have desired results in $B_r$ for $r<r_0$. For any $V'\in \IV(B_R)$ stationary with $F_{B_R}(V',\C)<\epsilon_0 R^n$ and a $\|V'\|$-measurable $E\subset \R^{n+k}$, we write $\tilde V'=\eta_{\lambda^{-1}\#}V'$, $\tilde E=\eta_{\lambda^{-1}}(E)$ for $\lambda=2\frac{R}{r_0}$. Since $\eta_{\lambda^{-1}\#}\C=\C$, we have obtain $\mathbf{F}_{B_{r_0/2}}(\tilde V',\C)\le \lambda^{-n} \mathbf{F}_{B_R}(V',\C)<\epsilon_0(r_0/2)^n$. We apply Theorem~\ref{relative.isoperimetric.varifold} to $\tilde V'$ and $\tilde E$ and scale back to get the desired result.
\end{proof}
\subsection{Local Harnack inequality on annuli}
Next, we would like to prove a Harnack inequality on annuli analogous to Theorem~\ref{LocalHarnack.balls}. This would give estimates for harmonic functions on annuli around an indecomposable point which provide estimates for the Green function. This theorem does not follow directly from Theorem~\ref{LocalHarnack.balls} by covering $A_{\beta_1,r}$ by balls. This is because if two Euclidean balls $B_1$, $B_2$ intersect, $\spt V\cap B_1$ may not intersect $ \spt V\cap B_2$. Thus the number of balls in a chain of balls that connect one point to another may depend on $r$.  Note also that we would need the extra assumption that any $\C\in \tan(V,x)$ satisfies $\HH^{n-1}(\sing\C)=0$ in order to apply the indecomposablility criterion. Recall that we write $A_{t,r}(x)=B_{(1+t)r}(x)\setminus \overline{B_{r/(1+t)}(x)}$.
\begin{theorem}[Local Harnack inequality on annuli]\label{LocalHarnack.annulus}
	Let $V\in \IV_n(\R^{n+k})$ be stationary  and $x\in\indec V$. Assume that any $\C\in \tan(V,x)$ satisfies $\HH^{n-1}(\sing\C)=0$. Then there exists $\beta_1=\beta_1(n,V,x)\in(0,1)$, $r_1=r_1(n,V,x)>0$, $\epsilon_1=\epsilon_1(n,V,x)\in(0,1)$ with the following significance. For any $r<r_1$, $V'\in \IV_n(A_{1,r}(x))$ stationary, $\mathbf{F}_{A_{1,r}(x)}(V,V')<\epsilon_1 r^n$, $\spt V'\cap A_{\beta_1/2, r}(x)\neq\emptyset$ and $0\le u\in W^{1,2}( A_{1,r}(x),V')$ satisfying $\Delta_{V'}u=0$ in $A_{1,r}(x)$ we have
	\[\sup _{\spt V\cap A_{\beta_1,r}(x)}u\le C\inf_{\spt V\cap A_{\beta_1,r}(x)}u\]
    and
    \[[u]_{\alpha,\spt V'\cap A_{\beta_1, r}(x)}\le Cr^{-\alpha}\|u\|_{L^\infty(A_{1,r}(x),\|V'\|)}\]
    for some $C=C(n,k,p,\Lambda,\beta_1)$ and $\alpha=\alpha(n,k,p,\Lambda,\beta_1)\in (0,1)$ where $\Lambda=\omega_n^{-1}r_0^{-n}\|V'\|(A_{1,r_0}(x))$.
\end{theorem}
Similar to the proof of Theorem~\ref{LocalHarnack.balls}, we need the following analogue of Lemma~\ref{volbd}.
\begin{lemma}\label{volbd.annulus}
    Let $V\in \IV(A_{1,r})$ be such that $\|V\|(A_{t,r})^{(n-1)/n}> 2\gamma\|\delta V\|(A_{t,r})$ for all $0<a<t<b\le 1$. Then $\|V\|(A_{b,r})\ge c(b-a)^nr^n$ where $c=c(n,k)=(8\gamma n)^{-n}$ and $\gamma$ is the constant from Lemma~\ref{isoperimetric}.
\end{lemma}
\begin{proof}
	By Definition~\ref{distributional.boundary}, we have
	\[  V\partial A_{t,r}=V\partial B_{\frac{r}{1+t}}-V\partial B_{(1+t)r}.\]
	By the coarea formula Theorem~\ref{coarea}, we have
	\[
	\begin{aligned}
		\frac{d}{dt}\|V\|(A_{t,r})=&\frac{d}{dt}\int_{\frac{r}{1+t}}^{(1+t)r} \int _{A_{1,r}}\frac{1}{|\nabla ^V|x||}d\|V\partial B_{\tau}\|d\tau\\
		= &r\int_{A_{1,r}} \frac{1}{|\nabla ^V|x||}d\|V\partial B_{(1+t)r}\| + \frac{r}{(1+t)^2}\int_{A_{1,r}} \frac{1}{|\nabla ^V|x||}d\|V\partial B_{\frac{r}{1+t}}\|\\
		\ge& r\|V\partial B_{(1+t)r}\|(A_{1,r}) +\frac{r}{(1+t)^2}\|V\partial B_{\frac{r}{1+t}}\|(A_{1,r})\\
		\ge &\frac{r}{4}\left(\|V\partial B_{(1+t)r}\|(A_{1,r}) +\|V\partial B_{\frac{r}{1+t}}\|(A_{1,r})\right)\\
        \ge&\frac{r}{4} \|V\partial A_{t,r}\|(A_{1,r}).
	\end{aligned}
	\]
	Thus using Theorem~\ref{isoperimetric}, we get
	\[\|V\|(A_{t,r})^{\frac{n-1}{n}}\le \|\delta V\|(A_{t,r})+\frac{4}{r}\frac{d}{dt}\|V\|(A_{t,r}).\]
	By our assumption,
	\[\frac{1}{2\gamma}\|V\|(A_{t,r})^{\frac{n-1}{n}}\le \frac{4}{r}\frac{d}{dt}\|V\|(A_{t,r}).\]
	Integrating the inequality, we get the conclusion.
\end{proof}
\begin{proof}[Proof of Theorem~\ref{LocalHarnack.annulus}]
    First we prove an analogue of Theorem~\ref{relative.isoperimetric.varifold} with $B_{tr}$ replaced by $A_{t,r}$ for any $0<t\le 1$ and $r>0$. We argue by contradiction and proceed exactly as Theorem~\ref{relative.isoperimetric.varifold} replacing Lemma~\ref{volbd} by Lemma~\ref{volbd.annulus}. We obtain $\C\in\tan(V,x)$ satisfying $\C\llcorner A_{1,1}=\C_1+\C_2$,  $\delta\C_i=0$ and $\|\C_i\|(A_{1,1})\ge c>0$ for $i=1,2$. Hence $\C\llcorner A_{1,1}$ is decomposable. By Lemma~\ref{reg.connected.mult1.equiv.indecomposable}, $\reg \C\cap A_{1,1}$ is disconnected. Since $\reg \C\cap A_{1,1}=\Sigma\times (\frac{1}{2},2)$ topologically where $\Sigma=\reg \C\cap \SSS^{n+k-1}$, we must have $\Sigma$ itself is disconnected. Thus $\reg \C=\Sigma\times (0,\infty)$ is disconnected which contradicts $\C$ is indecomposable and Lemma~\ref{reg.connected.mult1.equiv.indecomposable}.

    Once the analogue of Theorem~\ref{relative.isoperimetric.varifold} on annuli is established, one obtains a corresponding Neumann--Sobolev inequality (Theorem~\ref{Neumann.sobolev.inequality}) on annuli. The proof of the Harnack inequality then follows from the same argument as in Theorem~\ref{LocalHarnack.balls}.
\end{proof}

\begin{remark}
    If $V$ is the associated varifold of an area-minimizing boundary, then $r_1(n,V,x)=\infty$ for any $x\in \spt V$ for the same reason as Remark \ref{area.min.harnack}.
\end{remark}
Similar to Corollary~\ref{cone.r_0=infty}, we have the following.
\begin{corollary}\label{cone.r_1=infty}
    Let $\C\in \IV_n(\R^{n+k})$ be a stationary cone which is indecomposable at $0\in \spt \C$ and $\HH^{n-1}(\sing\C)=0$. Then Theorem~\ref{LocalHarnack.annulus} holds with $r_1(n,\C,0)=\infty$.
\end{corollary}
\section{Convergence of functions on converging varifolds}\label{convergence.of.functions}
We discuss $L^p$-convergence and uniform convergence of functions defined on the supports of converging varifolds. In this setting, we will prove the Rellich compactness theorem (Lemma~\ref{rellich.varying.varifolds}) and the Arzel\`a-Ascoli theorem (Lemma~\ref{arzela.ascoli}), corresponding to these respective modes of convergence. 
We write $\reg_\delta V$ for the set $\{x\in \reg V:d(x,\sing V)> \delta\}=\spt V\setminus(\sing V)_\delta$ for $\delta>0$.

\subsection{\texorpdfstring{$L^p$}{Lp} convergence}

\begin{defnlem}\label{Gromov.Hausdorff.Approximation}
	Let $V,\{V_i\}_{i=1}^\infty$ be stationary integral $n$-varifolds in $U$. Assume that $V_i\to V$ in $\IV_n(U)$ and that $V$ satisfies \eqref{mult1}.  Then for any $i\in \N$, there exist $\delta_i>0$ and a smooth diffeomorphism onto its image $\Phi_i:\reg_{\delta_i} V\to \reg V_i$ so that $\delta_i\to0$ as $i\to \infty$, $\sup_i\|\Phi_i\|_{C^k(\reg_{\delta_i} V)}\le C(k)$ for any $k\in\N$ and for any $\delta>0$ fixed, $\Phi_i|_{\reg_\delta V}\to \id|_{\reg_{\delta}V}$ smoothly.
\end{defnlem}
\begin{proof}
	By (\ref{mult1}), Allard's regularity theorem \cite[Theorem 24.2]{Sim} and elliptic regularity \cite{GT}, for any $x\in \reg V$, there is $r_x>0$  such that $\spt V_i\cap B_{r_x}(x)$ converges to $\spt V\cap B_{r_x}(x)$ smoothly as normal graphs. We take $\delta_i$ so that the graphing map $\Phi_i$ of $\spt V_i$ over $\spt V$ is well-defined on $\reg_{\delta_i} V$ and $\sup_i\|\Phi_i\|_{C^k(\reg_{\delta_i} V)}\le C(k)$ for any $k\in\N$. Since for every $x\in \reg V$, $\spt V_i\cap B_{r_x}(x)$ is a graph over $\spt V\cap B_{r_x}(x)$ for $i$ large enough, we can arrange $\delta_i\to0$. Since $\spt V_i\cap B_{r_x}(x)\to\spt V\cap B_{r_x}(x)$ smoothly for any $x\in \reg V$, we have $\Phi_i|_{\reg_\delta V}\to \id|_{\reg_{\delta}V}$ smoothly for $\delta>0$ fixed.
\end{proof}

We introduce the notion of $L^p$ convergence of functions on a sequence of converging multiplicity one stationary varifolds cf. \cite[Definition~1.2]{D}.

\begin{definition}\label{lp}
	Let $1 \le p < \infty$, $W \ssubset U$ and $V, \{V_i\}_{i=1}^\infty$ be stationary integral $n$-varifolds in $U$. Assume that $V_i \to V$ in $\IV_n(U)$ and that $V$ satisfies  $\HH^{n}(\sing V)=0$ and \eqref{mult1}.
    For each $i=1,2,\dots,$ let $u_i\in L^p(W,\|V_i\|)$ and $u\in L^p(W,\|V\|)$. We say that $u_i\to u$ in $L^p$ on $W$ as $V_i\to V$ if for any $\epsilon>0$ there exists $\delta>0$ so that $u_i|_{\Phi_i(\reg_\delta V\cap W)}\circ \Phi_i\to u|_{\reg_\delta V\cap W}$ in $L^p(W,\|V\|)$, \[\limsup_{i\to\infty}\|u_i|_{\spt V_i\cap W\setminus \Phi_i(\reg_\delta V\cap W)}\|_{L^p(W,\|V_i\|)}< \epsilon\] and \[\|u|_{\spt V\cap W\setminus \reg_\delta V}\|_{L^p(W,\|V\|)}< \epsilon.\]
\end{definition}

We have the following criterion for $L^p$-convergence.
\begin{lemma}\label{lp.convergence.criterion}
	Let $1 \le p < \infty$, $W \ssubset U$ and $V, \{V_i\}_{i=1}^\infty$ be stationary integral $n$-varifolds in $U$. Assume that $V_i \to V$ in $\IV_n(U)$ and that $V$ satisfies $\HH^{n}(\sing V)=0$ and \eqref{mult1}. If  for some $p'>p$ we have $u\in L^{p'}(W,\|V\|)$,  $u_i\in L^{p'}(W,\|V_i\|)$ with $\sup_i\|u_i\|_{L^{p'}(W,\|V_i\|)}<\infty$ and $u_i|_{\Phi_i(\reg_\delta V\cap W)}\circ \Phi_i\to u|_{\reg_\delta V\cap W}$ in $L^p(W,\|V\|)$ for any $\delta>0$ then $u_i\to u$ in $L^p$ on $W$ as $V_i\to V$.
\end{lemma}
\begin{proof}
	For any $\delta>0$, we define $S_\delta=\{x\in U:d(x,\sing V\cap \bar W)\le \delta\}$. By Lemma~\ref{Gromov.Hausdorff.Approximation}, we have $\Phi_i|_{\reg_\delta V}\to \id|_{\reg_{\delta}V}$ smoothly as $i\to\infty$. Hence we have 
	\[
		\|V_i\|(\reg_{2\delta}V_i\cap W\setminus\Phi_i(\reg_\delta V\cap W))\to0\text{ as }i\to\infty
	\]
	and $\spt V_i\cap W=  ( \reg_{2\delta} V_i\cap W) \cup (\spt V_i\cap S_{3\delta}\cap W)$ for $i$ sufficiently large. Thus 
	\[
	\spt V_i\cap W\setminus \Phi_i(\reg_\delta V\cap W)\subset (\reg_{2\delta}V_i\cap W\setminus\Phi_i(\reg_\delta V\cap W) )\cup (S_{3\delta}\cap W).
	\]
	By $\HH^{n}(\sing V)=0$, we have $\|V\|(S_{3\delta}\cap \bar W)\to0 $ as $\delta\to0$. Since $V_i\to V$ in $\IV_n(U)$, we have $\limsup\limits_{i\to\infty}\|V_i\|(S_{3\delta}\cap \bar W)\le \|V\|(S_{3\delta}\cap \bar W)$. We choose $\delta$ small enough such that 	\[\|u|_{\spt V\cap W\setminus \reg_\delta V}\|_{L^p(W,\|V\|)}\le \|V\|(S_{\delta}\cap \bar W)^{1-p/p'}\|u\|_{L^{p'}(W,\|V\|)}^{p/p'}<\epsilon\]
	and
	\[\begin{aligned}
		&\limsup_{i\to\infty}\|u_i|_{\spt V_i\cap W\setminus \Phi_i(\reg_\delta V\cap W)}\|_{L^p(W,\|V_i\|)}\\
		\le& \limsup_{i'\to\infty}(\|V_i\|(S_{3\delta}\cap \bar W)+\|V_i\|(\reg_{2\delta}V_i\cap W\setminus\Phi_i(\reg_\delta V\cap W)))^{1-p/p'}\|u_i\|_{L^{p'}(W,\|V_i\|)}^{p/p'}<\epsilon.
	\end{aligned} \]
	This verifies Definition~\ref{lp}.
\end{proof}

\begin{corollary}\label{div.X.converge}
		Let $1 \le p < \infty$ and $V, \{V_i\}_{i=1}^\infty$ be stationary integral $n$-varifolds in $U$. Assume that $V_i \to V$ in $\IV_n(U)$ and that $V$ satisfies $\HH^{n}(\sing V)=0$ and \eqref{mult1}. Then for any $W\ssubset U$ and $X\in C_c^2(U;\R^{n+k})$, we have $\Div_{V_i} X\to \Div_{V}X$ in $L^p$ on $W$ as $V_i\to V$.
\end{corollary}	
\begin{proof}
	By Lemma~\ref{Gromov.Hausdorff.Approximation}, we have $V_i\to V$ smoothly on $\reg_\delta V$ and hence $\Div_{V_i}X|_{\Phi_i(\reg_\delta V\cap W)}\circ \Phi_i\to \Div_{V}X|_{\reg_\delta V\cap W}$ in $L^p(W,\|V\|)$ for any $\delta>0$. Moreover, $\sup_i\|\Div_{V_i}X\|_{L^{p'}(W,\|V_i\|)}\le \sup_i\|V_i\|(W)^{1/p'}\|X\|_{C^2(\R^{n+k})}<\infty$. Thus  Lemma~\ref{lp.convergence.criterion} implies $\Div_{V_i}X\to \Div _VX$ in $L^p$ on $W$ as $V_i\to V$.
\end{proof}
Now we generalize Theorem~\ref{rellich.on.one.varifold} to converging multiplicity one varifolds.
\begin{lemma}[Rellich compactness for converging varifolds]\label{rellich.varying.varifolds}
	Let $1\le p<\infty$ and $V,\{V_i\}_{i=1}^\infty$ be stationary integral $n$-varifolds in $U$. Assume that $V_i \to V$ in $\IV_n(U)$ and that $V$ satisfies $\HH^{n}(\sing V)=0$ and \eqref{mult1}. Suppose $u_i\in W^{1,p}(V_i)$ are such that $\sup_i\|u_i\|_{W^{1,p}(V_i)}<\infty$. Then for any $W\ssubset U$  there exist $u\in L^{q}(W,\|V\|)$ and a subsequence $\{i'\}$ of $\{i\}$ so that $u_{i'}\to u$ in $L^q$ on $W$ as $V_{i'}\to V$ for $1\le q<\frac{np}{n-p}$.
\end{lemma}
\begin{proof}
	Let $1\le q<q'<\frac{np}{n-p}$ and $W\ssubset W_1\ssubset W_2\ssubset U$ be fixed. We define $S_\delta=\{x\in U:d(x,\sing V\cap \bar W_2)\le \delta\}$. We take $\zeta_\delta\in C_c^{\infty}(\R^{n+k})$ such that $\zeta_\delta=1$ on $\R^{n+k}\setminus S_{\delta}$ and $\zeta_\delta =0$ on $S_{\delta/2}$, $0\le\zeta_\delta\le 1$. We define 
	\begin{equation}\label{g_idelta}
		g_i^\delta=\zeta_{\delta} u_i|_{\Phi_i(\reg_{\delta/2} V\cap W_2)}\circ \Phi_i:\reg V\cap W_2\to \R.
	\end{equation}
    We have $g_i^\delta\in W^{1,p}_\loc(W_2\setminus \sing V,V)$. By construction, $g_i^\delta=0$ in $ W_2\cap (\sing V)_{\delta/2}$.  Lemma~\ref{sobolev.function.vanishing.on.singular.set} applied to $g_i^\delta$ on $V\llcorner W_2$ gives $g_i^\delta\in W^{1,p}(W_1,V)$. By (\ref{g_idelta}), Lemma~\ref{Gromov.Hausdorff.Approximation} and the $W^{1,p}$ bound on $u_i$, we have 
	\begin{equation}\label{w1p.bound.gdelta}
		\sup_i\|g_i^\delta\|_{W^{1,p}(W_1,V)}\le C(\delta).
	\end{equation} 
	Thus by Theorem~\ref{rellich.on.one.varifold}, there is a subsequence $\{i'\}$ such that 
	\begin{equation}\label{g.i.delta.to.g.delta}
		g_{i'}^\delta\to g^\delta \text{ in } L^{q'}(W,\|V\|)\text{ for any }\delta>0.
	\end{equation} By  $0\le \zeta_\delta \le 1$, we have for any $\delta>0$
	\[\begin{aligned}
		&\|g^\delta\|_{L^{q'}(W,\|V\|)}=\lim_{i\to\infty}\|g_i^\delta\|_{L^{q'}(W,\|V\|)}\\
        \le&\limsup_{i\to\infty}\|u_i|_{\Phi_i(\reg_{\delta/2} V\cap W_2)}\circ \Phi_i\|_{L^{q'}(W,\|V\|)} \\
		\le& \sup_i\|u_i\|_{L^{q'}(W,\|V_i\|)}\limsup_{i\to\infty}\|\Phi_i\|_{C^1(\reg_{\delta /2}V)}\le C
	\end{aligned}
	\]
	where in the last inequality we used Lemma~\ref{Gromov.Hausdorff.Approximation} and Theorem~\ref{Michael.Simon.Sobolev.inequality}. Note that the constant $C$ is independent of $\delta$. For $\delta_1>\delta_2$, we have $g^{\delta_1}_i=g^{\delta_2}_i$ on $\reg _{\delta _1}V\cap W$ and hence $g^{\delta_1}=g^{\delta_2}$ on $\reg _{\delta _1}V\cap W$. We define $u=g^\delta$ on $\reg_\delta V\cap W$ for any $\delta>0$. Then we have $\|u\|_{L^{q'}(W,\|V\|)}\le C$. By a diagonal argument, we can choose a subsequence $\{i'\}$ and $\delta_{i'}\to0$ so that $g_{i'}^{\delta_{i'}}\to u$ in $ L^{q'}(W,\|V\|)$. Since $\zeta_\delta=1$ on $\reg_\delta V$, we have $u_{i'}|_{\Phi_{i'}(\reg_\delta V\cap W)}\circ \Phi_{i'}\to u|_{\reg_\delta V\cap W}$ in $L^{q'}(W,\|V\|)$ for all $\delta>0$. Then Lemma~\ref{lp.convergence.criterion} implies $u_{i'}\to u$ in $L^q$ on $W$ as $V_{i'}\to V$.
\end{proof}
\subsection{Uniform convergence}
We introduce the concept of uniform convergence for converging stationary varifolds. Since $\spt V_i\to\spt V$ locally in Hausdorff distance by \cite[Chapter 4 Theorem 7.4]{simon2014introduction}, the concept of uniform convergence is the same as that of functions on metric spaces converging in Hausdorff distance.
\begin{definition}
	Let $V,\{V_i\}_{i=1}^\infty$ be stationary integral $n$-varifolds in $U$ such that $V_i\to V$ in $\IV_n(U)$. Let $u_i$ and $u$ be functions on $\spt V_i$ and $\spt V$ respectively. We say that $u_i\to u$ uniformly on compact subsets in $U$ if for any $\epsilon>0$ and $W\ssubset U$ compact, there are $\delta>0$ and $N\in\N$ such that for any $x\in \spt V\cap W$, $i>N$ and $x_i\in \spt V_i\cap W$ with $|x-x_i|<\delta$, we have $|u_i(x_i)- u(x)|<\epsilon$.

	A sequence of functions $u_i:\spt V_i\to \R$ is said to be locally equicontinuous if for any $\epsilon>0$ and $W\ssubset U$ compact there is $\delta>0$ such that for any $i\in\N$, $x,y\in \spt V_i\cap W$ with $|x-y|<\delta$, we have $|u_i(x)-u_i(y)|< \epsilon$.
\end{definition}
We have the Arzela-Ascoli theorem for equicontinuous family of functions on converging varifolds.
\begin{lemma}\label{arzela.ascoli}
	Let $V,\{V_i\}_{i=1}^\infty$ be stationary integral $n$-varifolds in $U$ such that $V_i\to V$ in $\IV_n(U)$ and $u_i:\spt V_i\to \R$ be a locally equicontinuous family of functions. Then there is a subsequence and a continuous function $u$ such that $u_i\to u$ uniformly on compact subsets.
\end{lemma}
\begin{proof}
	 Since $\spt V_i\to\spt V$ locally in Hausdorff distance by \cite[Chapter 4 Theorem 7.4]{simon2014introduction}, the lemma follows from the Arzela-Ascoli Theorem for metric spaces converging in Gromov-Hausdorff topology. See \cite[Lemma 45]{petersen2006riemannian}.
\end{proof}
The main result of this section is the following lemma showing that solutions to the Poisson equation converge under convergence of varifolds to a varifold with multiplicity one and codimension 2 singular sets.
\begin{lemma}\label{convergence.solution.poisson.equation}
	Let $V,\{V_i\}_{i=1}^\infty$ be stationary integral $n$-varifolds in $\R^{n+k}$ satisfying \eqref{euclidean.volume.growth}. Assume that $V_i \to V$ in $\IV_n(\R^{n+k})$ and that $V$ satisfies \eqref{Hn-2(sing)<infty} and \eqref{mult1}. Assume that either $u_i\in W_0^{1,2}( \Omega_i,V_i)$ for some exhaustion $\{\Omega_i\}_{i=1}^{\infty}$ of $\R^{n+k}$ by bounded open subsets or $u_i\in W_\loc^{1,2}(V_i)$ and $u_i\to0$ as $|x|\to\infty$ and $u_i$  satisfy $-\Delta_{V_i} u_i=f$ in the sense of (\ref{weaksln}) for $f\in C_c^{\infty}(\R^{n+k})$. Suppose $u_i\to u$ in $L^{q}$ on compact subsets as $V_i\to V$ for $1\le q<\frac{2n}{n-2}$. Then $u\in W^{1,2}_\loc(V)\cap L^\infty(\R^{n+k},\|V\|)$, $u\to0$ as $|x|\to\infty$ and $-\Delta_V u=f$ in the sense of (\ref{weaksln}). Moreover, for any $x\in \indec V$, there exists $r_1=r_1(n,V,x)>0$ such that $u_i\to u$ uniformly in $B_{r_1}(x)$ as $i\to\infty$. If in addition $V=\C$ is a cone which is indecomposable at 0, then $u_i\to u$ uniformly in $B_{r}$ for any $r>0$.
\end{lemma}
\begin{proof}We prove the case where $u_i\in W_\loc^{1,2}(V_i)$ and $u_i\to0$ as $|x|\to\infty$. The case where $u_i\in W_0^{1,2}(\Omega_i,V_i)$ is similar with minor modification.

	First we show that $u\in L^\infty(\R^{n+k},\|V\|)$. By Theorem~\ref{apri}, we have for any $R>0$, $\|u_i\|_{W^{1,2}( B_{4R},V_i)}\le C(R)$ and $\|u_i\|_{L^{\infty}(\R^{n+k},\|V_i\|)}\le C$. Since $u_i\to u$ in $L^{q}$ on compact subsets as $V_i\to V$ for $1\le q<\frac{2n}{n-2}$, we have $u\in L^{\infty}(\R^{n+k},\|V\|)$.

	Next, we show that $u\in W^{1,2}_\loc(V)$. 
	Let $R>0$ be arbitrary. We take $W=B_{R}$, $W_1=B_{2R}$, $W_2=B_{3R}$ and $U=B_{4R}$. By Lemma~\ref{cac} and Theorem~\ref{apri}, we have for any $R>0$, $\|u_i\|_{W^{1,2}( B_{4R},V_i)}\le C(R)$. Then the assumptions of Lemma~\ref{rellich.varying.varifolds} are satisfied. We show that the subsequential limit constructed in the proof of Lemma~\ref{rellich.varying.varifolds} is in $W^{1,2}(W,V)$. We adopt the same notations as in the proof of Lemma~\ref{rellich.varying.varifolds}. By (\ref{w1p.bound.gdelta}) and the weak compactness of $W^{1,2}(W_1,V)$, we can arrange, in addition to (\ref{g.i.delta.to.g.delta}), that $g_{i'}^\delta\to g^\delta$ weakly in $W^{1,2}(W_1,V)$. Since $\zeta_\delta=1$ on $\reg_\delta V\cap W_1$, we have $g_i^\delta=u_i|_{\Phi_i(\reg_{\delta/2} V\cap W_1)}\circ \Phi_i$ on $\reg_\delta V\cap W_1$. Then by weak convergence and Lemma~\ref{Gromov.Hausdorff.Approximation}, we have 
	\[\begin{aligned}
		\|g^\delta\|_{W^{1,2}(W_1\setminus (\sing V)_\delta,V)}\le& \limsup_{i\to\infty}\|g_i^\delta\|_{W^{1,2}(W_1\setminus (\sing V)_\delta,V)}\\
		\le& \sup_i\|u_i\|_{W^{1,2}(U,V_i)}\limsup_{i\to\infty}\|\Phi_i\|_{C^2(\reg_\delta V\cap U)}\le C.
	\end{aligned}\] 
	Note that the bound is independent of $\delta>0$.  Since $u_i\to u$ in $L^{q}$ on $W_1$ as $V_i\to V$, we have $ u=g^\delta$ on $\reg_\delta V\cap W_1$ for any $\delta>0$. Thus we have $u\in W^{1,2}_\loc(W_1\setminus\sing V,V)$ and $\|\nabla^Vu\|_{L^2(W_1,\|V\|)}\le C$. Since  $\|u\|_{L^{\infty}(W_1,\|V\|)}\le C$, we can apply Lemma~\ref{getting.sobolev.functions} to see that $u\in W^{1,2}(W,V)$. 
	
	Now we show that $u$ satisfies $-\Delta_Vu=f$. Since $u_i\in W^{1,2}_\loc(V_i)\cap L^\infty(\R^{n+k},\|V_i\|)$, by Proposition~\ref{integration.by.parts} and (\ref{weaksln}) we have
	\[\int u_i\Div_{V_i}(\nabla\phi) d\|V_i\|=-\int f\phi d\|V_i\|.\]
	Since $V_i\to V$, we have \[\lim_{i\to\infty}\int f\phi d\|V_i\|= \int f\phi d\|V\|.\]
	Since $u_i\to u$ in $L^{q}$ on compact subsets as $V_i\to V$ for $1\le q<\frac{2n}{n-2}$, and $\Div_{V_i}(\nabla\phi)\to \Div_{V}(\nabla\phi)$ in $L^{p}$ on compact subsets as $V_i\to V$ for $1\le p<\infty$ by Corollary~\ref{div.X.converge}, we have
	\[\lim_{i\to\infty}\int u_i\Div_{V_i}(\nabla\phi)  d\|V_i\|=\int u\Div_{V}(\nabla\phi) d\|V\|.\]
	Thus we get 
	\[\int u\Div_{V}(\nabla\phi) d\|V\|=-\int f\phi d\|V\|.\]
	Since $u\in W^{1,2}_\loc(V)\cap L^\infty(\R^{n+k},\|V\|)$, by Proposition~\ref{integration.by.parts} we have 
	\[\int \nabla ^V u\cdot \nabla^V\phi d\|V\|=\int f\phi d\|V\|.\]

    Finally we show uniform convergence. By Theorem~\ref{LocalHarnack.balls}, we have $\|u_i\|_{C^\alpha(\spt V_i\cap B_r(x))}\le C$ for any $r<\beta_0(n,V,x) r_0(n,V,x)$. Lemma~\ref{arzela.ascoli} implies that there is a subsequence $u_{i'}$ converges to $u$ uniformly. If $V=\C$ is an indecomposable cone, then $r$ can be taken arbitrarily since $r_0(n,\C,0)=\infty$ by Corollary~\ref{cone.r_0=infty}.
\end{proof}

\section{Green functions on stationary integral varifolds}\label{greens.function.stationary.varifolds}
In this section, we construct the Dirichlet Green function and the global Green function on a stationary varifold with Euclidean volume growth following Gr\"uter-Widman \cite{gruter1982green}. We note that in Proposition~\ref{DiriGreen} and Theorem~\ref{globalGreenvarifolds}, the pole can be any point $x\in \spt V$, whether regular or singular. We assume $n\ge3$.
    
\begin{proposition}\label{DiriGreen}
	Let $V\in \IV_n(\R^{n+k})$ be stationary and satisfy \eqref{euclidean.volume.growth}. Then for any $\Omega\subset \R^{n+k}$ bounded open, we have the following.
	\begin{enumerate}
		\item[(i)]  There exists a function $G_\Omega:(\spt V\cap \Omega)\times (\spt V\cap \Omega)\to [0,\infty]$
		satisfying: for any $x\in \spt V\cap \Omega$, $r>0$ we have 
		\begin{equation}\label{w1p.dirichlet}
			G_\Omega(x,\cdot)\in W^{1,2}( \Omega\setminus B_{r}(x),V)\cap W^{1,p}_0(\Omega,V)\quad \forall p\in [1,\frac{n}{n-1})
		\end{equation}
		and  $\Delta_{V,y} G_\Omega(x,y)=-\delta_x(y)$ in the sense that for any $\phi\in C_c^\infty(\Omega)$, 
		\begin{equation}\label{delta.dirichlet}
			\int\nabla _y^V G_\Omega(x,y) \cdot \nabla^V\phi(y)d\|V\|(y)=\phi(x).
		\end{equation}
	\end{enumerate}
	Moreover, the function $G_\Omega$ has the following properties.
	\begin{enumerate}
		\item[(ii)] For any $x\in \spt V\cap \Omega$, $\|G_\Omega(x,\cdot)\|_{L_*^{n/(n-2)}(\Omega,\|V\|)}\le C(n)$.
		
		\item[(iii)] For any $x\in \spt V\cap \Omega$, $\|\nabla^V G_\Omega(x,\cdot)\|_{L_*^{n/(n-1)}(\Omega,\|V\|)}\le C(n,\Lambda)$.
		
		\item[(iv)] For  any $x\in \spt V\cap \Omega$ and $\|V\|$-a.e. $y\in\spt V\cap \Omega$ or any $y\in (\reg V\cup \indec V)\cap \Omega$, we have $G_\Omega(x,y)\le C(n,\Lambda)|x-y|^{2-n}$.
		
		\item[(v)] For $\|V\|$-a.e. $x,y\in \spt V\cap\Omega$ or any $x,y\in (\reg V\cup \indec V)\cap \Omega$ we have $G_\Omega(x,y)=G_\Omega(y,x)$.
		
		\item[(vi)]Let $u\in W_{0}^{1,2}(\Omega,V)$ be the solution to $-\Delta_V u=f$ where $f\in C_c^\infty (\Omega)$. Then for $\|V\|$-a.e. $x\in \spt V\cap \Omega$ or any $x\in(\reg V\cup \indec V)\cap \Omega$, we have 
		\[u(x)=\int G_\Omega(x,y)f(y)d\|V\|(y).\]
	\end{enumerate}
\end{proposition}

\begin{proof}
    The proof of the proposition follows the argument of \cite[Theorem 1.1]{gruter1982green}, with several straightforward modifications. First, volume bounds for balls are obtained using \eqref{euclidean.volume.growth}. Second, the Euclidean Sobolev inequality is replaced by Theorem~\ref{Michael.Simon.Sobolev.inequality}, and the classical mean value inequality is replaced by Theorem~\ref{mean.value.inequality}. Finally, whenever a pointwise estimate or pointwise convergence is used in \cite{gruter1982green}, the corresponding statement holds for $\|V\|$-a.e. $x \in \spt V\cap\Omega$ by standard measure-theoretic arguments, and also holds for every $x \in (\reg V \cup \indec V)\cap \Omega$ by elliptic regularity \cite{GT} together with Theorem~\ref{LocalHarnack.balls}.

	We now sketch the main steps in the proof. Let $x\in \spt V\cap \Omega$ be fixed. By the Riesz representation theorem, for any $\rho>0$ there exists a unique $G^\rho=G^\rho(x,\cdot)\in W_0^{1,2}(\Omega,V)$ such that for any $\phi\in W_0^{1,2}(\Omega,V)$
	\begin{equation}
		\int \nabla G^\rho\cdot \nabla\phi d\|V\|=\fint_{B_\rho (x)}\phi d\|V\|:=\frac{1}{\|V\|(B_\rho(x))}\int_{B_\rho(x)}\phi d\|V\|.\label{approximate.solution}
	\end{equation}
    We call such function the approximate Green function of $V$ on $\Omega$. Following exactly the same argument in \cite[Theorem 1.1]{gruter1982green} we get, among other things, $\|G^\rho\|_{L_*^{n/(n-2)}(\Omega,\|V\|)}\le C(n)$ and $\|\nabla^V G^\rho\|_{L_*^{n/(n-1)}(\Omega,\|V\|)}\le C(n,\Lambda)$.   Hence for any $1\le p<\frac{n}{n-1}$, $\|G^\rho\|_{W^{1,p}(\Omega,V)}\le C(n,p,\Lambda,\Omega)$.
	By Theorem~\ref{rellich.on.one.varifold} and weak compactness of $W^{1,p}(\Omega,V)$, there exists $\rho_j\to0$ so that \begin{equation}\label{Grho.to.GOmega}
		G^{\rho_j}\to G_\Omega=G_\Omega(x,\cdot)
	\end{equation} strongly in $L^q(\Omega,\|V\|)$ for $1\le q<\frac{n}{n-2}$ and weakly in $W^{1,p}_0(\Omega,V)$.
    
	One can verify (i)-(iv) using the same argument as \cite[Theorem 1.1]{gruter1982green}.
	
	To show (v), let $\rho_i,\sigma_j\to0$ be such that $G^{\rho_i}(x,\cdot)\to G_\Omega(x,\cdot)$ and $G^{\sigma_j}(y,\cdot)\to G_\Omega(y,\cdot)$ strongly in $L^1(\Omega,\|V\|)$, $\|V\|$-a.e. in $\spt V\cap \Omega$ and uniformly around regular and indecomposable points by elliptic regularity \cite{GT} and Theorem~\ref{LocalHarnack.balls}. We take $G^{\rho_i}(x,\cdot)$ and $G^{\sigma_j}(y,\cdot)$ as test functions in (\ref{approximate.solution}) for $G^{\sigma_j}(y,\cdot)$ and $G^{\rho_i}(x,\cdot)$ respectively and get
	\[
		\fint_{B_{\rho_i}(x)}G^{\sigma_j}(y,\cdot)d\|V\|=\fint_{B_{\sigma_j}(y)}G^{\rho_i}(x,\cdot)d\|V\|.
	\]
	We consider the limit of both sides by first letting $i\to\infty$ and then $j\to\infty$.
	
	For the left hand side, letting $i\to\infty$ and using Lebesgue-Besicovitch differentiation theorem \cite[Theorem 1.32]{gariepy2015measure}, we have for any $y\in \spt V\cap\Omega$ and $\|V\|$-a.e. $x\in \spt V\cap \Omega$, 
	\[
		\lim_{i\to\infty}\fint_{B_{\rho_i}(x)}G^{\sigma_j}(y,\cdot)d\|V\|=G^{\sigma_j}(y,x).
	\]  Since for any $y\in \spt V\cap \Omega$, $G^{\sigma_j}(y,\cdot)\to G_\Omega(y,\cdot)$ $\|V\|$-a.e. in $\spt V\cap \Omega$, we have for $\|V\|$-a.e. $x\in \spt V\cap \Omega$
	\[
		\lim_{j\to\infty}\lim_{i\to\infty}\fint_{B_{\rho_i}(x)}G^{\sigma_j}(y,\cdot)d\|V\|=\lim_{j\to\infty}G^{\sigma_j}(y,x)=G_\Omega(y,x).
	\]
	
	For the right hand side, letting $i\to\infty$ since $G^{\rho_i}(x,\cdot)\to G_\Omega(x,\cdot)$ in $L^1(\Omega,\|V\|)$, we get for any $x,y\in \spt V\cap \Omega$ 
	\[
		\lim_{i\to\infty}\fint_{B_{\sigma_j}(y)}G^{\rho_i}(x,\cdot)d\|V\|=\fint_{B_{\sigma_j}(y)}G_\Omega(x,\cdot)d\|V\|.
	\]
	By the Lebesgue-Besicovitch differentiation theorem \cite[Theorem 1.32]{gariepy2015measure} again, for any $x\in \spt V\cap \Omega$ and $\|V\|$-a.e. $y\in\spt V\cap \Omega$, 
	\[
		\lim_{j\to\infty}\lim_{i\to\infty}\fint_{B_{\sigma_j}(y)}G^{\rho_i}(x,\cdot)d\|V\|=\lim_{j\to\infty}\fint_{B_{\sigma_j}(y)}G_\Omega(x,\cdot)d\|V\|=G_\Omega(x,y).
	\]
	Thus we get $G_\Omega(x,y)=G_\Omega(y,x)$ for $\|V\|\times\|V\|$-a.e. $(x,y)\in (\spt V\cap \Omega)\times (\spt V\cap \Omega)$.
	The same identity is true if $x,y\in (\reg V\cup \indec V)\cap \Omega$ since $G_\Omega(x,\cdot),G^\rho(x,\cdot)$ are continuous around $y$ and $G_\Omega(y,\cdot),G^\rho(y,\cdot)$ are continuous around $x$ by elliptic regularity \cite{GT} and Theorem~\ref{LocalHarnack.balls}.
	This finishes (v).
	
	To show (vi), since $u\in W_0^{1,2}(\Omega,V)$ solves $-\Delta_Vu=f$ and $G^\rho(x,\cdot)\in W_0^{1,2}(\Omega,V)$, we have
	\[\int\nabla^V u(y)\cdot\nabla^V G^\rho(x,y)d\|V\|(y)=\int G^\rho(x,y)f(y)d\|V\|(y).\]
	By (\ref{approximate.solution}) we have
	\[\fint_{B_\rho(x)}ud\|V\|=\int G^\rho(x,y)f(y)d\|V\|(y).\]
	Taking $\rho\to0$, and using Lebesgue-Besicovitch differentiation theorem \cite[Theorem 1.32]{gariepy2015measure} and $G^\rho(x.\cdot)\to G_\Omega(x,\cdot)$ in $L^1(\Omega,\|V\|)$, we get for $\|V\|$-a.e. $x\in \spt V\cap \Omega$
	\[
	u(x)=\int G_\Omega(x,y)f(y)d\|V\|(y).
	\]
	If $x\in (\reg V\cup \indec V)\cap \Omega$, then the same identity holds at $x$ since $u$ is continuous around $x$ by elliptic regularity \cite{GT} and Theorem~\ref{LocalHarnack.balls}.
\end{proof}
\begin{theorem}\label{globalGreenvarifolds} 
	Let $V\in \IV_n(\R^{n+k})$ be stationary and satisfy \eqref{euclidean.volume.growth}.
	\begin{enumerate}
		\item[(i)]  There exists a function $G:\spt V\times \spt V\to [0,\infty]$
		satisfying: for any $x\in \spt V$, $r>0$ we have 
		\begin{align}
			&G(x,\cdot)\in W_\loc^{1,2}( \R^{n+k}\setminus B_{r}(x),V)\cap W_\loc^{1,p}(V)\quad \forall p\in [1,\frac{n}{n-1})\label{w1p.global}\\
			&G(x,y)\to0\quad\text{as }|y|\to\infty \label{decay.global}
		\end{align}
		and  $\Delta_y G(x,y)=-\delta_x(y)$ in the sense that for any $\phi\in C_c^\infty(\R^{n+k})$, 
		\begin{equation}\label{delta.global}
			\int\nabla _y^V G(x,y) \cdot \nabla^V\phi(y)d\|V\|(y)=\phi(x).
		\end{equation}
	\end{enumerate}
	Moreover, the function $G$ has the following properties.
	\begin{enumerate}
		\item[(ii)] For any $x\in \spt V$, $\Omega\subset\R^{n+k}$ bounded open subset and $1\le p<\frac{n}{n-1}$, $\|G(x,\cdot)\|_{W^{1,p}(\Omega,V)}\le C(n,p,\Omega,\Lambda)$.
		
		\item[(iii)] For any $x\in \spt V$ and $\|V\|$-a.e. $y\in\spt V$ or any $y\in \reg V\cup \indec V$, we have
		\begin{equation}\label{rough.upper.bound}
			G(x,y)\le C(n,\Lambda)|x-y|^{2-n}.
		\end{equation}
		
		\item[(iv)] For $\|V\|$-a.e. $x,y\in \spt V$ or any $x,y\in \reg V\cup \indec V$ we have
		\begin{equation}
			\label{symmetry} G(x,y)=G(y,x).
		\end{equation}
		\item[(v)]  Let $u\in W_\loc^{1,2}(V)$ be the solution to $-\Delta_V u=f$ and $u\to0$ as $|x|\to\infty$ where $f\in C_c^\infty (\R^{n+k})$. 
		Then for $\|V\|$-a.e. $x\in \spt V$ or any $x\in \reg V\cup \indec V$
		\begin{equation}\label{representation} 
			u(x)=\int G(x,y)f(y)d\|V\|(y).
		\end{equation}
	\end{enumerate}
\end{theorem}

\begin{proof}
	Let $x\in \spt V$ be fixed. Let $\{\Omega_i\}_{i=1}^\infty$ with $\Omega_i\subset \Omega_{i+1}$ be an exhaustion of $\R^{n+k}$ by bounded open subsets. By Proposition~\ref{DiriGreen} (i), for each $i\in \N$, there exists a Dirichlet Green function $G_{\Omega_i}(x,\cdot)$ on $ V\llcorner \Omega_i$. By Proposition~\ref{DiriGreen} (ii), (iii) we have for any $R>0$ $\|G_{\Omega_i}(x,\cdot)\|_{W^{1,p}(B_R(x),V)}\le C(n,p,R,\Lambda)$ for $1\le p<\frac{n}{n-1}$ and any $i\in \N$. By Theorem~\ref{rellich.on.one.varifold}, weak compactness of the Banach space $W^{1,p}(B_R(x),V)$ and a diagonal argument, there exists a subsequence still denoted by $\{i\}$ and $G(x,\cdot)\in W^{1,p}_\loc(V)$, so that 
	\begin{equation}\label{G_Omega.to.G}
		G_{\Omega_i}(x,\cdot) \to G(x,\cdot)
	\end{equation} strongly in  $L^q(B_R(x),\|V\|)$, $\|V\|$-a.e. in $\spt V\cap B_R(x)$ and weakly in $W^{1,p}(B_R(x),V)$ and $W^{1,2}(B_R(x)\setminus B_r(x),V)$ for $1\le p<\frac{n}{n-1}$, $1\le q<\frac{n}{n-2}$ and every $R>r>0$. Moreover, if $y\in \reg V\cup \indec V\setminus\{x\}$, then $G_{\Omega_i}(x,\cdot)\to G(x,\cdot)$ uniformly around $y$ as $i\to\infty$ by elliptic regularity \cite{GT} and Theorem~\ref{LocalHarnack.balls}. By the weak convergence and that $G_{\Omega_i}$ satisfies Proposition~\ref{DiriGreen} (i), we see that (\ref{w1p.global}), (\ref{delta.global}) and (ii) are satisfied.
	
	Since $G_{\Omega_i}(x,\cdot)\ge0$ we have $G(x,\cdot)\ge0$.
	
	By Proposition~\ref{DiriGreen} (iv) and \eqref{G_Omega.to.G}, we have $G(x,y)\le C(n,\Lambda)|x-y|^{2-n}$. Hence (\ref{decay.global}) and (iii) are satisfied.
	
	(iv) Since for each $i$, $G_{\Omega_i}(x,y)=G_{\Omega_i}(y,x)$ for $\|V\|$-a.e. $x,y\in \spt V\cap \Omega_i$ or any $x,y\in (\reg V\cup \indec V)\cap \Omega_i$, we have $G(x,y)=G(y,x)$ for $\|V\|$-a.e. $x,y\in \spt V$ or any $x,y\in \reg V\cup \indec V$ by the $\|V\|$-a.e. convergence and uniform convergence around regular and indecomposable points.

	(v) Let $f\in C_c^\infty (\R^{n+k})$ be fixed and $u_i\in W_0^{1,2}(V)$ be the solution to $-\Delta_V u_i=f$  on $V\llcorner \Omega_i$. By Lemma~\ref{cac} and Theorem~\ref{apri}, $\|u_i\|_{W^{1,2}( B_R(x),V)}\le C(n,k,\Lambda,R,f)$. By Theorem~\ref{rellich.on.one.varifold}, weak compactness of the Banach spaces $W^{1,p}( B_R(x),V)$ and a diagonal argument, there exists a subsequence still denoted by $\{i\}$ and $u\in W^{1,2}_\loc(V)$ such that $u_i\to u$ strongly in $L^2(B_R(x),\|V\|)$, $\|V\|$-a.e. in $\spt V\cap B_R(x)$ and weakly in $W^{1,2}( B_R(x),V)$. Moreover, if $x\in \reg V\cup \indec V$, then $u_i\to u$ smoothly around $x$ by elliptic regularity \cite{GT} and Theorem~\ref{LocalHarnack.balls}. Thus by weak convergence, $u$ solves $-\Delta_V u=f$. For $y\in \spt V$ and $|x|\ge 2R:=2d(0,\spt f)$, we have $|x-y|\ge |x|-R\ge \frac{1}{2}|x|$. By Proposition~\ref{DiriGreen} (iv) and (vi), for $i\in \N$ large and $|x|\ge 2R$, we have
	\begin{equation}\label{proof.decay.for.poisson.equgation}
		|u_i(x)|\le C(n,\Lambda)\int |x-y|^{2-n}|f(y)|d\|V\|(y)\le C(n,\Lambda)(\frac{1}{2})^{2-n}\|f\|_{L^1(\|V\|)}|x|^{2-n}.
	\end{equation}
	Since $u_i\to u$ for $\|V\|$-a.e. $x\in \spt V$, we have $u\to0$ as $|x|\to\infty$. Note that solutions decaying at infinity are unique by Theorem~\ref{apri}. Thus $u$ is \textit{the} solution to $-\Delta_V u=f$ and $u\to0$ as $|x|\to\infty$. On the other hand, by construction $G_{\Omega_i}(x,\cdot)\to G(x,\cdot)$ in $L^q(B_R(x),\|V\|)$, for any $x\in \spt V$ and 
	$R>0$. We take $R$ large such that $\spt f\subset B_R(x)$. Then we have $\int G_{\Omega_i}(x,y)f(y)d\|V\|(y)\to \int G(x,y)f(y)d\|V\|(y)$ for every $x\in \spt V$. Combining the convergence of $u_i$, we get $u(x)=\int G(x,y)f(y)d\|V\|(y)$ for $\|V\|$-a.e. $x\in \spt V$ and for any $x\in \reg V\cup \indec V$.
\end{proof}
By (\ref{proof.decay.for.poisson.equgation}), we have the following.
\begin{corollary}\label{decay.for.poisson.equgation}
	Let $V\in \IV_n(\R^{n+k})$ be stationary and satisfy \eqref{euclidean.volume.growth}. Let $u\in W_\loc^{1,2}(V)$ be the solution to $-\Delta_V u=f$ and $u\to0$ as $|x|\to\infty$ where $f\in C_c^\infty (\R^{n+k})$.
	Then for $\|V\|$-a.e. $x\in \spt V$ or any $x\in \reg V\cup \indec V$ such that $|x|\ge 2d(0,\spt f)$
	\[|u(x)|\le C(n,\Lambda,\|f\|_{L^1(\|V\|)})|x|^{2-n}.\]
\end{corollary}

As an example, we can explicitly compute the Green function constructed as above on a stationary varifold which is also a cone.
\begin{lemma}\label{cone}
	Let $\C\in \IV_n(\R^{n+k})$ be a stationary varifold which is also a cone. Then the Green function constructed in Theorem~\ref{globalGreenvarifolds} with the pole at $0$ is 
	\[G(0,y)=\frac{1}{n(n-2)\omega_n\Theta(\C,0)}|y|^{2-n}.\]
\end{lemma}
\begin{proof}
	To simplify notations, we write $c=\frac{1}{n(n-2)\omega_n\Theta(\C,0)}$ and $r=|y|$. It suffices to show that the approximate Green function $G^\rho$ of $\C$ on $B_R$ in the proof of Theorem~\ref{DiriGreen} is given by 
	\begin{equation}\label{Grho.on.cone}
		G^\rho (y)=\left\{\begin{aligned}
			&\frac{1}{2}c\rho^{-n}(n\rho^2-(n-2)r^2)-cR^{2-n} &\quad &r\le \rho\\
			&cr^{2-n}-cR^{2-n} &\quad &\rho<r<R.
		\end{aligned}\right.
	\end{equation}
	If this is done, then using (\ref{Grho.on.cone}) and the fact from (\ref{Grho.to.GOmega}) that $G^\rho\to G_{B_R}$ in $W^{1,p}_0( B_R,\C)$ as $\rho\to0$, the Dirichlet Green function of $\C$ on $B_R$ is
	\begin{equation}\label{GBR.on.cone}
		G_{B_R}=cr^{2-n}-cR^{2-n}.
	\end{equation}
	Moreover, using (\ref{GBR.on.cone}) and the fact from (\ref{G_Omega.to.G}) that $G_{B_R}\to G$ in $W^{1,p}_\loc(\C)\cap W^{1,2}_\loc (\R^{n+k}\setminus B_r,\C)$ as $R\to\infty$, we have $G=cr^{2-n}$ which is the desired result.
	
	To see (\ref{Grho.on.cone}). We note that $G^{\rho}$ is the restriction of a Lipschitz function on $\R^{n+k}$ which vanishes on $\partial B_R$. Thus $G^\rho \in W^{1,2}_0( B_R,\C)$. We compute \[\nabla G^{\rho}=\left\{\begin{aligned}
		&(2-n)c\rho^{-n}y&\quad r<\rho\\
		&(2-n)cr^{-n}y&\quad r>\rho,
	\end{aligned}\right.
	\text{ and }\Div_\C(\nabla G^{\rho})= \left\{\begin{aligned}
		&n(2-n)c\rho^{-n}&\quad r<\rho\\
		&0&\quad r>\rho.
	\end{aligned}\right.\]
	For any $\phi\in W_0^{1,2} ( B_R,V)$, using Proposition~\ref{integration.by.parts} and the fact that $\|\C\|(B_\rho)=\omega_n \Theta(\C,0)\rho^n$, we have 
	\[\int \nabla G^\rho\cdot \nabla^\C \phi d\|\C\|=-\int \Div_\C (\nabla G^\rho)\phi d\|\C\|=\frac{1}{\omega_n \Theta(\C,0)\rho^n}\int_{B_\rho}\phi d\|\C\|=\fint _{B_\rho}\phi d\|\C\|.\]
	Thus we see that $G^\rho$ solves (\ref{approximate.solution}). Since such function is unique in $W_0^{1,2}(B_R,\C)$ by Theorem~\ref{apri}, we have $G^{\rho}$ is the approximate Green function of $\C$ on $B_R$.
\end{proof}

\section{Convergence of Green functions}\label{convergence.of.green.function}
We prove the following convergence theorem of Green functions in this section.
\begin{theorem}\label{conv1}
	Let $V,\{V_i\}_{i=1}^\infty$ be stationary integral $n$-varifolds in $\R^{n+k}$ satisfying \eqref{euclidean.volume.growth}. Assume that $V_i \to V$ in $\IV_n(\R^{n+k})$ and that $V$ satisfies \eqref{Hn-2(sing)<infty} and \eqref{mult1}. Let $\{\Omega_i\}_{i=1}^\infty$ be either an exhaustion of $\R^{n+k}$ by bounded open subsets, or $\Omega_i = \R^{n+k}$ for all $i$. Let $G$ be the Green function of $V$ and $G_i$ be the Green function of $V_i$ on $\Omega_i$. Assume that $\spt V_i\ni x_i\to x_0\in \reg V\cup \indec V$. 
    
    Then we have $G_i(x_i,\cdot)\to G(x_0,\cdot)$ in $L^q$ on compact subsets as $V_i\to V$ for $1\le q<\frac{n}{n-2}$. Moreover, the convergence is uniform in a neighborhood of any $y\in \indec V\setminus \{x_0\}$ and in any small annulus centered at $x_0$.
    
    If $V=\C$ is a cone and $x_0=0$, then we know further that $G_i(x_i,\cdot)\to G(0,\cdot)$ uniformly in any annulus centered at $0$.
\end{theorem}
\begin{proof}
    The proof for the case where $\{\Omega_i\}_{i=1}^\infty$ is an exhaustion of bounded open sets and where $\Omega_i=\R^{n+k}$ are similar. The necessary modification are to use Theorem~\ref{DiriGreen} instead of Theorem~\ref{globalGreenvarifolds}. We present here the proof for $\Omega_i=\R^{n+k}$.
    
	Let $f\in C_c^\infty(\R^{n+k})$ and $u_i\in W^{1,2}_\loc(V_i)$ be such that $-\Delta_{V_i} u_i=f$ and $u_i\to0$ as $|x|\to\infty$. Then by Lemma~\ref{cac} and Theorem~\ref{apri}, we have for any $R>0$, $\|u_i\|_{W^{1,2}(B_R,V_i)}\le C(n,k,\Lambda,R,f)$ and $\|u_i\|_{L^{\infty}(B_R,\|V_i\|)}\le C(n,k,\Lambda,R,f)$.  Thus by Lemma~\ref{rellich.varying.varifolds}, Lemma~\ref{arzela.ascoli} and a diagonal argument, there is a subsequence such that $u_{i'}\to u$ in $L^{q}$ on compact subsets as $V_{i'}\to V$ for $1\le q<\frac{2n}{n-2}$.
    By Lemma~\ref{convergence.solution.poisson.equation}, we have $u\in W^{1,2}_\loc(V)$ satisfies $-\Delta _V u=f$ in the sense of (\ref{weaksln}) and $u_{i'}\to u$ uniformly in $B_r(x_0)$ for $r>0$ sufficiently small. By Corollary~\ref{decay.for.poisson.equgation} we have $|u_i(x)|\le C(n,\Lambda,\|f\|_{L^1(\|V_i\|)})|x|^{2-n}$ for $|x|\ge 2d(0,\spt f)$. Thus by $\|V\|$-a.e. convergence, we have $u\to 0$ as $|x|\to\infty$.
	
	On the other hand, for any $R>0$, by Theorem~\ref{globalGreenvarifolds} (ii), $\|G_i(x_i,\cdot)\|_{W^{1,p}(B_R,V_i)}\le C(n,k,\Lambda,R)$ for $1\le p<\frac{n}{n-1}$. By Theorem~\ref{LocalHarnack.balls}, for any $y\in (\reg V\cup\indec V) \setminus\{x_0\}$ there is $r_2=r_2(n,V,y)<r_0(n,V,y)\beta_0(n,V,y)$ such that \[\|G_i(x_i,\cdot)\|_{C^\alpha(\spt V_i\cap B_{r_2}(y))}\le C(n,k,\Lambda,\beta_0)\|G_i(x_i,\cdot)\|_{L^\infty( B_{r_2}(y),\|V_i\|)}.\] Moreover, by Theorem~\ref{LocalHarnack.annulus}, there is $r_3=r_3(n,V,x_0)=\beta_1(n,V,x_0)r_1(n,V,x_0)$ such that for any $r<r_3$, \[\|G_i(x_i,\cdot)\|_{C^\alpha(\spt V_i\cap A_{\beta_1,r}(x_0))}\le C(n,k,\Lambda,\beta_1)r^{-\alpha}\|G_i(x_i,\cdot)\|_{L^\infty( A_{\beta_1,r}(x_0),\|V_i\|)}.\] Thus by Lemma~\ref{rellich.varying.varifolds}, Lemma~\ref{arzela.ascoli} and a diagonal argument, there is a subsequence $\{i'\}$ and $g_{x_0}\in L^q_\loc(\R^{n+k},\|V\|)$ such that $G_{i'}(x_{i'},\cdot)\to g_{x_0}(\cdot)$ in $L^q$ on compact subsets as $V_{i'}\to V$ for $1\le q<\frac{n}{n-2}$ and $G_{i'}(x_{i'},\cdot)\to g_{x_0}(\cdot)$ uniformly in $B_{r_2}(y)$ and $A_{\beta_1,r}(x_0)$. 
	Hence $\int G_{i'}(x_{i'},y)f(y)d\|V_{i'}\|(y)\to \int g_{x_0}(y)f(y)d\|V\|(y)$. By Theorem~\ref{globalGreenvarifolds} (v), since $x_0\in \reg V\cup\indec V$, we have $u_{i'}(x_{i'})=\int G_{i'}(x_{i'},y)f(y)d\|V_{i'}\|(y)$ and $u(x_0)=\int G(x_0,y)f(y)d\|V\|(y)$ by Theorem~\ref{LocalHarnack.balls}. Thus \[\begin{aligned}
		&\int G(x_0,y)f(y)d\|V\|(y)=u(x_0)=\lim_{i'\to\infty} u_{i'}(x_{i'})\\
		=&\lim_{i'\to\infty} \int G_{i'}(x_{i'},y)f(y)d\|V_{i'}\|(y)=\int g_{x_0}(y)f(y)d\|V\|(y).
	\end{aligned}\]Since the above equality is true for all $f\in C_c^\infty(\R^{n+k})$, we have $g_{x_0}(\cdot)=G(x_0,\cdot)$ $\|V\|$-a.e. on $ \spt V$ and pointwise in $B_{r_2}(y)$ for $y\in (\reg V\cup\indec V) \setminus\{x_0\}$ and in $A_{\beta_1,r}(x_0)$ for $r<r_3$.  This shows that $G_{i'}(x_{i'},\cdot)\to G(x_0,\cdot)$ in the desired sense.
	
	One can repeat the argument in previous paragraphs to show that any subsequence of $\{G_i(x_i,\cdot)\}$ has a further subsequence converging to $G(x_0,\cdot)$. By the standard subsequence characterization of convergence, the original sequence $G_i(x_i,\cdot)$ converges to $G(x_0,\cdot)$ in the desired sense.
	
    The result when $V=\C$ is a cone follows from the above argument and Corollary~\ref{cone.r_0=infty} and \ref{cone.r_1=infty}.
\end{proof}
\section{Behavior of the Green function with the pole at an indecomposable point}\label{behavior.green.function.singular.point}
\subsection{Asymptotics of the Green function near the indecomposable pole or infinity}
Using the convergence theorem established in the previous section, we analyze the behavior of the Green function on a stationary varifold near the indecomposable pole and at infinity.
\begin{theorem}\label{asymp1}
    Let $V\in \IV_n(\R^{n+k})$ be stationary and satisfy \eqref{euclidean.volume.growth} and $\Omega\subset \R^{n+k}$ be a bounded open set. Denote by $G_\Omega(x,y)$ the Dirichlet Green function of $V$ on $\Omega$ and by $G(x,y)$ the Green function of $V$.
    If $x\in \reg V$ or $x\in \indec  V$ such that any $\C\in \tan(V,x)$ satisfies \eqref{Hn-2(sing)<infty}, then we have
    \[\lim_{y\to x}G_\Omega(x,y)|y-x|^{n-2}=\lim_{y\to x}G(x,y)|y-x|^{n-2}=\frac{1}{n(n-2)\omega_n\Theta(V,x)}.
	\]
	If $V$ is indecomposable at $\infty$ and any $\C\in \tan(V,\infty)$ satisfies \eqref{Hn-2(sing)<infty}, then for any $x\in\spt V$
	\[\lim_{|y|\to \infty}G(x,y)|y-x|^{n-2}=\frac{1}{n(n-2)\omega_n\Theta(V,\infty)}.\]
\end{theorem}
\begin{proof}The proof of the three limits are similar with minor modifications. We only compute $\lim\limits_{y\to x}G(x,y)|y-x|^{n-2}$.
	Let $y_i\in \spt V\setminus \{x\}$ be such that $y_i\to x$ as $i\to\infty$. We set $\lambda_i=|y_i-x|$. Consider $V_i=(\eta_{x,\lambda_i^{-1}})_\# V$ and $\tilde y_i=\frac{y_i-x}{|y_i-x|}$. Let $G_i$ be the Green function on $V_i$. Then $G_i(0,\tilde y_i)=G(x,y_i)|y_i-x|^{n-2}$. Since $\lambda_i\to0$, by \cite[Theorem 42.7]{Sim} there exist a subsequence $\{i^\prime\}$, $\C\in \tan(V,x)$ and $y\in \spt \C\cap \SSS^{n+k-1}$ so that $V_{i^\prime}\to \C$ in $\IV(\R^{n+k})$ and $\tilde y_{i'}\to y$. By Lemma~\ref{reg.connected.mult1.equiv.indecomposable}, $\C$ satisfies \eqref{mult1}. By Theorem~\ref{conv1} (convergence on any annulus) and Lemma~\ref{cone} we have \[G(x,y_{i'})|y_{i'}-x|^{n-2}=G_{i'}(0,\tilde y_{i'})\to G(0,y)=\frac{1}{n(n-2)\omega_n\Theta(\C,0)}=\frac{1}{n(n-2)\omega_n\Theta(V,x)}.\] Since every sequence of $G(x,y_i)|y_i-x|^{n-2}$ has a subsequence converging to $\frac{1}{n(n-2)\omega_n\Theta(V,x)}$, we establish the limit.
\end{proof}
\subsection{Global upper and lower bounds of the Green function}
The above asymptotics gives lower and upper bounds of the Green function near an indecomposable pole or at infinity. Next, we investigate global upper and lower bounds of the Green function.

We have a sharp upper bound of the Green function with the pole at an indecomposable point. We need the following preliminary lemmas.
For $x,y\in \R^{n+k}$, we write the Euclidean Green function on $\R^n$ as
\[\Gamma(x,y)=\frac{1}{n(n-2)\omega_n}|x-y|^{2-n}.\]
\begin{lemma}\label{super}
	Let $V\in \IV_n(U)$ be stationary. Then for any $x\in \R^{n+k}$, $\phi\in C_c^\infty(U)$ with $\phi\ge0$ and $x\notin \spt \phi$, we have
	\[\int \nabla^V\Gamma(x,\cdot)\cdot \nabla^V\phi d\|V\|=\frac{1}{\omega_n}\int \phi (y)\frac{|(y-x)^\perp|^2}{|y-x|^{n+2}}d\|V\|(y)\ge0. \]
	 Moreover, if $y\in\reg V$, then we have
	 \[\Delta_{V,y} \Gamma(x,y)=-\frac{1}{\omega_n}\frac{|(y-x)^\perp|^2}{|y-x|^{n+2}}\]
	where $(y-x)^{\perp}$ is the projection of $y-x$ to $T_yV^\perp$.
\end{lemma}
\begin{proof}
	Without loss of generality, we take $x=0$ and denote $r=|y|$. We compute
	\begin{equation}\label{div.r^2-n}
		\begin{aligned}
			\Div _V(\nabla r^{2-n})=&(2-n)\Div _V(r^{-n}x)=n(2-n)r^{-n}+(2-n)(-n)r^{-n}|\nabla ^Vr|^2\\
			=&n(2-n)r^{-n}|\nabla^\perp r|^2\\
			=&n(2-n)\frac{|y^\perp|^2}{|y|^{n+2}}.
		\end{aligned}
	\end{equation}
	For any $\phi\in C_c^\infty(U)$ with $\phi\ge0$ and $0\notin \spt \phi$, the function $r^{2-n}$ is smooth in a neighborhood of $\spt\phi$. Since $V$ is stationary, we have
	\[0=\int \Div _V(\phi\nabla r^{2-n})d\|V\|=\int \nabla^V r^{2-n}\cdot\nabla^V\phi d\|V\|+\int \phi \Div_V (\nabla r^{2-n})d\|V\|.\]
	Thus we get
	\[\int \nabla^V r^{2-n}\cdot\nabla^V\phi d\|V\|=n(n-2)\int \phi \frac{|y^\perp|^2}{|y|^{n+2}}d\|V\|\ge0.\]
	If $y\in \reg V$, we can take $\{e_i\}_{i=1}^n$ to be an orthonormal frame of $TV$ near $y$. Stationarity implies the mean curvature $\sum_{i=1}^n(\nabla_{e_i}e_i)^\perp=0$ near regular points. Hence
	\[\begin{aligned}
		\Div _V(\nabla r^{2-n})=&\Div _V(\nabla^V r^{2-n})+\Div _V(\nabla^\perp r^{2-n})=\Delta_Vr^{2-n}+\sum_{i=1}^n \nabla_{e_i}\nabla^\perp r^{2-n}\cdot e_i\\
		=&\Delta_Vr^{2-n}-\sum_{i=1}^n \nabla^\perp r^{2-n}\cdot \nabla_{e_i}e_i\\
		=&\Delta_Vr^{2-n}.
	\end{aligned}\]
	Combining with (\ref{div.r^2-n}), we get the second equation.
\end{proof}
We also need the following lemma about unique continuation of stationary varifolds.
\begin{lemma}\label{unique.continuation}
    Let $V_1,V_2\in\IV_n(U)$ be stationary such that $\HH^{n-1}(\sing V_2)=0$. If $\reg V_1$ is connected and $\reg V_1\cap W=\reg V_2\cap W$ for some $W\subset V$ open, then $\reg V_1\subset \spt V_2$.
\end{lemma}
\begin{proof}
    Let $x\in \reg V_1\cap W$ and $y\in \reg V_1\setminus \sing V_2$. Since $\HH^{n-1}(\sing V_2)=0$ and $\reg V_1$ is connected, there is a continuous path $\sigma :[0,1]\to \reg V_1\setminus \sing V_2$ such that $\sigma(0)=x$, $\sigma(1)=y$. Let 
    \[E=\{t\in[0,1]:\exists r>0, \reg V_1\cap B_r(\sigma(t))\subset \reg V_2\cap B_r(\sigma(t))\}.\] Since $0\in E$, $E$ is nonempty. $E$ is open by definition. Suppose $t_i\in E$ and $t_i\to t$. Then $\sigma(t_i)\in \reg V_2$ and $z=\lim\limits_{i\to\infty} \sigma(t_i)=\sigma(t)\in \spt V_2$. Since $\sigma$ takes value in $\reg V_1\setminus\sing V_2$, we have $z\in \reg V_1\cap \reg V_2$. Since $T_{\sigma(t_i)}V_1=T_{\sigma(t_i)}V_2$, we have $T_{z}V_1=T_zV_2$. There is $r>0$ such that $\reg V_1\cap B_r(z)$ is a graph of $u_1$ over $A_1\subset T_z V_1$ open and $\reg V_2\cap B_r(z)$ is a graph of $u_2$ over $A_2\subset T_z V_2$ open, where $0\in A_1\cap A_2$ and $z=u_1(0)=u_2(0)$.  Since $t_i\in E$ and $\sigma(t_i)\to z$, there is an open set in $A_1\cap A_2$ on which $u_1=u_2$. By unique continuation of elliptic systems \cite{Protter60}, $u_1=u_2$ on $A_1\cap A_2$. Thus $t\in E$ and $E$ is closed. By connectedness of $[0,1]$, $E=[0,1]$. Thus $y\in \reg V_2$. Then we have $\reg V_1\setminus \sing V_2\subset \reg V_2$ and hence $\reg V_1\subset \spt V_2$.
\end{proof}
    \begin{theorem}
        \label{sharp.upper.bound.singular.points}Let $V\in \IV_n(\R^{n+k})$ be stationary and satisfy \eqref{euclidean.volume.growth}. Denote by $G(x,y)$ the Green function of $V$. Let $x\in \reg V$ or $x\in \indec  V$ such that any $\C\in \tan(V,x)$ satisfies \eqref{Hn-2(sing)<infty}. Then we have for $\|V\|$-a.e. $y\in \spt V$ or $y\in \reg V\cup \indec V$
	\begin{equation}\label{upperglobal}
	    G(x,y)\le\frac{1}{n(n-2)\omega_n\Theta(V,x)}|x-y|^{2-n}.
	\end{equation}
	If $\HH^{n-1}(\sing V)=0$ and equality in \eqref{upperglobal} holds at $x\neq y\in \reg V\cup \indec V$, then $V=\eta_{x,1\#}\C_0+V_1$ where $\C_0$ is a stationary cone and $y-x\in \spt \C_0$ and $V_1\in \IV_n(\R^{n+k})$ is stationary with $x\notin \spt V_1$.
    \end{theorem}
\begin{proof}
    Let $x$ be as in the assumption and $\Omega\subset\R^{n+k}$ be bounded open. By Theorem~\ref{asymp1}, for $\epsilon>0$ there exists $r(\epsilon)>0$ satisfying $\lim_{\epsilon\to0}r(\epsilon)=0$ so that for any $|x-y|<r(\epsilon)$, $G_\Omega(x,y)<\frac{1}{\Theta(\|V\|,x)}\Gamma_{\epsilon,x}(y)$  where \[\Gamma_{\epsilon,x}(y)=\Gamma(x,y)+\epsilon|x-y|^{2-n}.\] We define $w_\epsilon(z)=\frac{1}{\Theta(\|V\|,x)}\Gamma_{\epsilon,x}(y)-G_\Omega(x,z)$. In view of Lemma~\ref{super}, we have $\Delta_V w_\epsilon\le 0$ in $\Omega\setminus B_r(x)$ for $r<r(\epsilon)$. By Theorem~\ref{apri}, we have $w_\epsilon\ge 0$ in $\Omega\setminus B_r(x)$. Letting $\epsilon\to0$, we get \begin{equation}\label{upper}
        G_\Omega(x,y)\le\frac{1}{n(n-2)\omega_n\Theta(V,x)}|x-y|^{2-n}.
    \end{equation}
    Let $\{\Omega_i\}_{i=1}^\infty$ with $\Omega_i\subset \Omega_{i+1}$ be an exhaustion of $\R^{n+k}$ by bounded open subsets.
    Since $G_{\Omega_i}$ satisfies \eqref{upper} for each $i$, and $G_{\Omega_i}\to G$ $\|V\|$-a.e. and uniformly near regular and indecomposable points, we see that $G$ satisfies \eqref{upperglobal}.

	We now discuss the rigidity case for (\ref{upperglobal}). Suppose $x\neq y\in \reg V\cup \indec V$ is such that the equality in \eqref{upperglobal} holds. We set \[w(z)=\frac{1}{\Theta(V,x)}\Gamma(x,z)-G(x,z).\] Then we have $w\ge 0$ and $w(y)=0$. By Lemma~\ref{super}, we have $\Delta_V w\le 0$ in $\R^{n+k}\setminus\{x\}$. Since $y\in \reg V\cup \indec V$, we apply Theorem~\ref{LocalHarnack.balls} to get $w\equiv 0$ in a neighborhood of $y$ in $\spt V$. Let $M=\{z\in \reg V:w(z)=0\}$. By Lemma~\ref{super} 
    \[ -\frac{1}{\omega_n}\frac{|(z-x)^\perp|^2}{|z-x|^{n+2}}=\Delta_{V,z} \Gamma(x,z)=\Theta (V,x)\Delta_{V,z}G(x,z)=0\quad \text{in }M.\]
	Thus $(z-x)^\perp=0$ for all $z\in M$. Integrating this equation, we see that $M$ is locally conical. Apply Lemma~\ref{unique.continuation} to each connected component of $\C(M)=\{x+r(z-x):z\in M\}$, we get $\C(M)\subset \spt V$. Let $\C_0\in \tan(V,x)$. Then $\C(M)\subset \eta_{x,1}(\spt \C_0)$. Since $\C_0$ satisfies \eqref{Hn-2(sing)<infty}, by Lemma~\ref{reg.connected.mult1.equiv.indecomposable}, $\reg\C_0$ is connected and $\C_0$ has multiplicity one. Thus $\C_0=|\reg \C_0|$. By Lemma~\ref{unique.continuation}, $\eta_{x,1}(\reg \C_0) \subset \spt V$. Thus $V=\eta_{x,1\#}\C_0+V_1$ for $V_1\in \IV_n(\R^{n+k})$. Since $V,\eta_{x,1\#}\C_0$ are stationary, we have $V_1$ is stationary. We must have $x\notin \spt V_1$ since otherwise, for any $\C\in \tan(V,x)$, we have $\C=\C_0+\C_1$ where $\C_1\in \tan(V_1,x)$. This contradicts the indecomposability of $\C$.
\end{proof}

In general, one does not have global lower bound for the Green function if the varifold is decomposable at some point or at infinity, see Example \ref{2Planesn-2}. When every point in $\spt V$ as well as infinity is indecomposable, we get a global lower bound for the Green function. We note that the constant for the lower bound is not explicit.

\begin{theorem}\label{global.lower.bound}
    Let $V\in \IV_n(\R^{n+k})$ be stationary and satisfy \eqref{euclidean.volume.growth} and $x\in \indec V$. Denote by $G(x,y)$ the Green function of $V$. Suppose $\spt V=\reg V\cup\indec V$ and $V$ is indecomposable at $\infty$ and that any $\C\in \tan(V,x)\cup  \tan(V,\infty)$ satisfies \eqref{Hn-2(sing)<infty}. Then there is $c=c(n,V,x)>0$ such that for any $y\in \spt V$
	\[ G(x,y)\ge c(n,V,x)|x-y|^{2-n}.\]
\end{theorem}
\begin{proof}
	Assume for contradiction that no such $c(n,V,x)>0$ exists. Then for each $i\in \N$, there exist $y_i\in \spt V$ so that $G(x,y_i)<i^{-1}|x-y_i|^{2-n}$. By Theorem~\ref{asymp1}, there are $R>r>0$ such that $r\le |y_i-x|\le R$ for $i$ large. By compactness of $\spt V\cap \bar B_R(x)\setminus B_r(x)$, there exists a subsequence $\{i'\}$ so that $ y _{i'}\to y\in \spt V\cap \bar B_R(x)\setminus B_r(x)$. Since $\spt V=\reg V\cup\indec V$, by Theorem~\ref{LocalHarnack.balls}, $G(x,\cdot)$ is continuous around $y$ and $G(x,y)=\lim_{i'\to\infty} G(x,y_{i'})=0$. By Theorem~\ref{LocalHarnack.balls}, $G(x,\cdot)\equiv 0$ in a neighborhood of $y$. Since $V$ is indecomposable at $\infty$, we have $\spt V$ is connected since otherwise the tangent cone at infinity is decomposable. Since  $\{y\in\spt V:G(x,y)= 0\}$ is nonempty, open and closed by  Theorem~\ref{LocalHarnack.balls}, we have $G(x,\cdot)\equiv 0$, a contradiction to Theorem~\ref{asymp1}.
\end{proof}
\begin{remark}\label{area.min.c=c(n)}
    As mentioned in Remark \ref{area.min.bdy.indecomposable}, if $V$ is the associated varifold of an area-minimizing boundary, then $V$ is indecomposable by \cite[Theorem 1]{bombieri1972harnack}. In particular,  every point in $\spt V$ as well as infinity is indecomposable. Moreover, the singular set of $V$ has codimension 7 which was discussed in Remark \ref{known.result.sing.v}.  Thus Theorem~\ref{conv1}, \ref{asymp1}, \ref{sharp.upper.bound.singular.points}, \ref{global.lower.bound} hold for area-minimizing boundaries. In fact, in this case, the constant $c$ in Theorem~\ref{global.lower.bound} does not depend on $V$ or $x$ since one can apply the same argument to extract a sequence $(V_i,x_i)$ violating it and repeat the same limiting argument to obtain a contradiction.
\end{remark}
As an application we have the following estimates for entire solutions to the minimal surface equation.
\begin{corollary}
	Let $u:\R^n\to \R$  be an entire solution to the minimal surface equation 
	\[(1+|\nabla u|^2)\Delta u-\sum_{i,j=1}^nu_{x_ix_j}u_{x_i}u_{x_j}=0.\]
	Then there exists $C=C(n,u)$ such that for $|x|\ge 1$
	\[|\nabla u(x)|\le C|x|^{n-2}\left(1+\frac{|u(x)-u(0)|}{|x|}\right)^{n-2}.\]
\end{corollary}

\begin{proof}
	Let $\Sigma=\graph u$. Translating $u$ vertically, we assume $u(0)=0$ and hence $0\in \Sigma$. We consider the function $v(x)=(1+|\nabla u|^2)^{-1/2}>0$. By \cite[2.8]{simon1987asymptotic} we have $\Delta_\Sigma v\le 0$.  By \cite[\S 3]{simon1983survey}, $\Sigma$ is an area-minimizing boundary. By Theorem~\ref{DiriGreen}, the Dirichlet Green function of $\Sigma$ on $B_R$ satisfies $G_{B_R}(0,\cdot)|_{\partial B_1\cap \Sigma}\le C(n)$. Let $m=\inf_{\partial B_1\cap\Sigma}v$. Then $mC(n)^{-1} G_{B_R}(0,\cdot)|_{\partial B_1\cap \Sigma}\le m$. Since $mC(n)^{-1}G_{B_R}(0,\cdot)\in W^{1,2}_0 ( B_R,V)$ is harmonic in $B_R\setminus B_1$, by Theorem~\ref{apri} we have $v\ge mC(n)^{-1} G_{B_R}$ in $\Sigma\cap B_{R}\setminus B_1$. Letting $R\to \infty$, we have $G_{B_R}(0,\cdot)\to G(0,\cdot)$ by (\ref{G_Omega.to.G}) and thus $v(x)\ge mC(n)^{-1} G(0,x)$ for $x\in\Sigma\setminus B_1$. By Theorem~\ref{global.lower.bound}, Remark \ref{area.min.c=c(n)}, we have
	\[v(x)\ge c(u(x)^2+|x|^2)^{-(n-2)/2},\]i.e.
	\[1+|\nabla u(x)|^2\le C(|x|^2+u(x)^2)^{n-2}.\]
	Thus we have for $|x|\ge 1$
	\[|\nabla u(x)|\le C|x|^{n-2}\left(1+\frac{|u(x)|}{|x|}\right)^{n-2}.\]
\end{proof}
\begin{remark}
	This estimate is stronger compared to a similar estimate of \cite[(6.8)]{bombieri1972harnack} since it is not known whether $|u(x)|\le c (1+|x|^{(n-2)/2})$ which is a long standing conjecture cf. \cite{bombieri1972harnack,Yau}. As mentioned in the final remark of \cite{bombieri1972harnack}, one would need the bound $\sup_{B_r}|\nabla u|\le C(1+\frac{1}{r}\sup_{B_r}|u|)$ to conclude polynomial growth of $u$ and our estimate here is not enough to deduce such estimate.
\end{remark}
\section{Green functions on decomposable varifolds}\label{remark.examples}

We discuss some results of the Green function on decomposable varifolds. The following uniqueness of the Green function with the pole at a regular point is useful.
\begin{lemma}\label{uniqueness.Green.regular.points}
	Let $\Omega\subset\R^{n+k}$ be a bounded open set and $V\in \IV_n(\R^{n+k})$ be stationary and satisfy \eqref{euclidean.volume.growth}.
	Then we have
	\begin{enumerate}
		\item [(i)]
		For any $x\in \reg V\cap \Omega$, $G_\Omega(x,\cdot)$ is the unique function satisfying (\ref{w1p.dirichlet}), (\ref{delta.dirichlet}). 
		\item[(ii)]
		For any $x\in \reg V$, $G(x,\cdot)$ is the unique function satisfying (\ref{w1p.global}), (\ref{decay.global}) and (\ref{delta.global}).
	\end{enumerate}
\end{lemma}
\begin{proof}
	(i) Suppose $G_{1,\Omega},G_{2,\Omega}$ are two functions satisfying (\ref{w1p.dirichlet}) and (\ref{delta.dirichlet}) for $x\in \reg V\cap \Omega$. By (\ref{w1p.dirichlet}), (\ref{delta.dirichlet}), $G_{1,\Omega}-G_{2,\Omega}\in W^{1,p}_0(\Omega,V)$ is harmonic in the sense of (\ref{weaksln}). Since $x\in \reg V$, by Weyl's lemma on smooth domains \cite{burch1978dini}, $G_{1,\Omega}-G_{2,\Omega}$ is smooth in $\spt V\cap B_r(x)$ for $r$ small. Combining with (\ref{w1p.dirichlet}), we have $G_{1,\Omega}-G_{2,\Omega}\in W^{1,2}_0(\Omega,V)$. By Theorem~\ref{apri} we have $G_{1,\Omega}-G_{2,\Omega}=0$.
	
	(ii) The uniqueness for the global Green function is similar. Suppose $G_1,G_2$ are two functions satisfying (\ref{w1p.global}), (\ref{decay.global}) and (\ref{delta.global}) for $x\in \reg V$. By (\ref{w1p.global}) and (\ref{delta.global}), $G_1-G_2\in W^{1,p}_\loc(V)$ is harmonic in the sense of (\ref{weaksln}). Since $x\in \reg V$, by Weyl's lemma on smooth domains \cite{burch1978dini}, $G_1-G_2$ is smooth in $\spt V\cap B_r(x)$ for $r$ small. Combining with (\ref{w1p.global}), we have $G_1-G_2\in W^{1,2}_\loc(V)$. Applying Theorem~\ref{apri} to $G_1-G_2$ on $V\llcorner B_R(x)$, using (\ref{decay.global}) and letting $R\to \infty$, we get $G_1-G_2=0$.
\end{proof}
\begin{corollary}\label{Green.V_1+V_2}
    Let $V_1,V_2\in \IV_n(\R^{n+k})$ be stationary and satisfy \eqref{euclidean.volume.growth}. Let $G_1,G_2$ be the Green function of $V_1,V_2$ and $G$ be the Green function of $V=V_1+V_2$. Assume that for any $W\ssubset\R^{n+k}$, $\HH^{n-2}(\spt V_1\cap \spt V_2\cap W)<\infty$, then for any $x\in \reg V_1\setminus\spt V_2$
    \begin{equation}\label{def.G.V_1+V_2}
        G(x,y)=\left\{\begin{aligned}
		&G_1(x,y)&\quad& y\in \spt V_1\setminus \spt V_2\\
		&0 &\quad &y\in \spt V_2\setminus \spt V_1.
	\end{aligned}\right.
    \end{equation} and for any $x\in \reg V_2\setminus\spt V_1$
    \[G(x,y)=\left\{\begin{aligned}
		&G_2(x,y)&\quad& y\in \spt V_2\setminus \spt V_1\\
		&0 &\quad &y\in \spt V_1\setminus \spt V_2.
	\end{aligned}\right.\]
\end{corollary}
\begin{proof}
    We only prove the case where $x\in \reg V_1\setminus \spt V_2$. In view of Lemma~\ref{uniqueness.Green.regular.points}, it suffices to show that $G(x,\cdot)$ as defined in \eqref{def.G.V_1+V_2} satisfies (\ref{w1p.global}), (\ref{decay.global}) and (\ref{delta.global}). For any $R>0$, let $S=\spt V_1\cap \spt V_2\cap B_R$. We have $\HH^{n-2}(S)<\infty$. By Theorem~\ref{globalGreenvarifolds}, $G(x,\cdot)\in W^{1,p}_\loc(B_R\setminus S,V)\cap  L^\infty( (S)_\delta\cap B_R,\|V\|)$ for $\delta<\frac{1}{2}d(x,S)$ and $\nabla^V G(x,\cdot)\in L^p(B_R, \|V\|)$ for $1\le p< \frac{n}{n-1}$. By Lemma~\ref{getting.sobolev.functions}, we have $G(x,\cdot)\in W^{1,p}(B_{R/2},V)$. Thus $G(x,\cdot)\in W^{1,p}_\loc (V)$. Similarly $G(x,\cdot)\in W^{1,2}_\loc (\R^{n+k}\setminus B_r(x),V)$ for all $r>0$. This verifies (\ref{w1p.global}). The other two properties (\ref{decay.global}) and (\ref{delta.global}) are obvious.
\end{proof}
\subsection{Discontinuity and breakdown of lower bounds for the Green function}
The following example shows that the Green function on a general stationary integral $n$-varifold can be discontinuous at decomposable points and is not bounded below by the Euclidean Green function. 

\begin{example}\label{2Planesn-2}
	Let $V=P_1+P_2$ be the union of two distinct $n$-planes $P_1,P_2$ in $\R^{n+k}$ such that $\dim P_1\cap P_2\le n-2$ and $G$ be its Green function. Since $P_1\cup P_2$ is a minimal cone, we can apply Lemma~\ref{cone} to see that for any $x\in P_1\cap P_2$
	\[G(x,y)=\frac{1}{2n(n-2)\omega_n}|x-y|^{2-n}.\]By Corollary~\ref{Green.V_1+V_2}, we have
	\[G(x,y)=\left\{\begin{aligned}
		&\Gamma(x,y)&\quad& x,y\in P_1\setminus P_2\text{ or }x,y\in P_2\setminus P_1\\
		&0 &\quad &\text{otherwise.}
	\end{aligned}\right.\]  
	
	We see that for $x\in P_1\cup P_2\setminus P_1\cap P_2$, $G(x,\cdot)$ is not continuous at $\sing V=P_1\cap P_2$ and $G(0,x)\neq G(x,0)$ for any $x\in P_1\cup P_2\setminus P_1\cap P_2$. As such, the restriction on \eqref{symmetry} and \eqref{representation} in Theorem~\ref{globalGreenvarifolds} being true $\|V\|$-a.e. and on $\reg V\cup \indec V$ is necessary. We also see that $G(x,y)=0$ for $x\in P_1\setminus P_2$ and $y\in P_2\setminus P_1$. Thus $G(x,\cdot)$ does not converge to $G(0,\cdot)$ in any possible sense as $ P_1\cup P_2\setminus P_1\cap P_2\ni x\to0$ and we don't have the asymptotics at infinity and the lower bound in Theorem~\ref{asymp1} and \ref{global.lower.bound}. The requirement that $x$ or $\infty$ is indecomposable in Theorem~\ref{conv1}, \ref{asymp1} and \ref{global.lower.bound} is necessary.
\end{example}

\subsection{Propagation of the Green function across a codimension-one singular set}
The following two examples shows that if the singular set is large i.e. codimension 1, then the Green function can ``propagate" from one component to another.
\begin{example}\label{2planesn-1}
	Let $V=P_1+P_2$ be the union of two distinct $n$-planes $P_1,P_2$ in $\R^{n+k}$ such that $\dim P_1\cap P_2=n-1$ and $G$ be its Green function. 
	As in Example \ref{2Planesn-2}, for any $x\in P_1\cap P_2$
	\[G(x,y)=\frac{1}{2n(n-2)\omega_n}|x-y|^{2-n}.\]
	Since $\dim P_1\cap P_2=n-1$, $P_1\cap P_2$ divides $P_1$ and $P_2$ into four components $P_{1}^+,P_1^-,P_2^+,P_2^-$. We identify $P_i^\pm$ with $\R^n_+=\{(x_1,\dots,x_n)\in \R^n:x_n>0\}$ in the following way. Let $\{e_1\dots,e_{n-1}\}$ be an orthonormal basis of $P_1\cap P_2$ and $\nu_i^\pm$ be the inward conormal of $P_1\cap P_2$ in $P_i^\pm$. We take the basis $\{e_1,\dots,e_{n-1},\nu_i^\pm\}$ on $P_i^\pm$. Then we can identify $P_i^\pm$ with $\R^n_+$. For any $x=(x_1,\dots,x_n)\in \R^n_+$, we write $x^*=(x_1,\dots,x_{n-1},-x_n)$.
	
	We claim that if $x\in P_1^+$ we have
	\[G(x,y)=\left\{\begin{aligned}
		&\Gamma(x,y)-\frac{1}{2} \Gamma (x^*,y)&\quad& y\in P_1^+\\
		&\frac{1}{2}\Gamma(x^*,y)&\quad&y\in P_1^-\cup P_2^+\cup P_2^- .
	\end{aligned}\right.\]
	where we view $x\in \R^n_+$ under the identification $P_1^+\cong \R^n_+$ and the formula on each component is written under the identification $P_i^\pm\cong \R^n_+$.
	The cases where $x$ is in other components can be similarly defined.
	
	To show the claim, we verify (\ref{w1p.global}), (\ref{decay.global}) and (\ref{delta.global}). Since $G(x,\cdot)$ is Lipschitz around $\sing V=P_1\cap P_2$, we see that $G(x,\cdot)$ satisfies (\ref{w1p.global}). It is clear that (\ref{decay.global}) is satisfied.
	We verify (\ref{delta.global}) as follows. We write $\nu$ as the inward conormal of $\R^n_+$. Note that if $y\in \R^{n-1}\times \{0\}$ then $\nu\cdot\nabla_y\Gamma(x,y)=-\nu\cdot\nabla_y \Gamma(x^*,y)$. We have
	\[\begin{aligned}
		&\int \nabla^V G(x,\cdot)\cdot\nabla^V \phi d\|V\|\\
		=&\phi(x)-\int \phi\nu\cdot \nabla \Gamma(x,\cdot)d\HH^{n-1}+\frac{1}{2}\int \phi\nu\cdot \nabla \Gamma(x^*,\cdot)  d\HH^{n-1}-\frac{3}{2}\int \phi\nu\cdot\nabla \Gamma(x^*,\cdot)  d\HH^{n-1}\\
		=&\phi(x).
	\end{aligned}\]
	In contrast to Example \ref{2Planesn-2}, $G(x,\cdot)$ is continuous for any $x\in \spt V$ and has a universal lower bound.
\end{example}

The following example compute the Green function on a stationary varifolds with triple junction singularities where similar phenomena occurs.
\begin{example}\label{triple.junction}
	Let $V=H_1+H_2+H_3$ be as in Lemma~\ref{triple.junction.indecomposable} and $G$ be its Green function. 
	We apply Lemma~\ref{cone} to see that for any $x\in H_1\cap H_2\cap H_3$
	\[G(x,y)=\frac{2}{3n(n-2)\omega_n}|x-y|^{2-n}.\]
	As in Example \ref{2planesn-1}, we identify $H_i$ with $ \R^n_+$ for $i=1,2,3$.
	We claim that if $x\in H_1\setminus (H_1\cap H_2\cap H_3)$ then we have
	\[G(x,y)=\left\{\begin{aligned}
		&\Gamma(x,y)-\frac{1}{3} \Gamma (x^*,y)&\quad& y\in H_1\\
		&\frac{2}{3}\Gamma(x^*,y)&\quad&\text{on }  y\in H_2\cup  H_3 
	\end{aligned}\right.\]
	Again, we verify (\ref{w1p.global}), (\ref{decay.global}) and (\ref{delta.global}). Since $G(x,\cdot)$ is Lipschitz around $\sing V=H_1\cap H_2\cap H_3$, we see that $G(x,\cdot)$ satisfies (\ref{w1p.global}). It is clear that (\ref{decay.global}) is satisfied.
	We verify (\ref{delta.global}) as follows. Using the same argument and notation as in Example \ref{2planesn-1}, we have
	\[\begin{aligned}
		&\int \nabla^V G(x,\cdot)\cdot\nabla^V \phi d\|V\|\\
		=&\phi(x)-\int \phi\nu\cdot \nabla \Gamma(x,\cdot)d\HH^{n-1}+\frac{1}{3}\int \phi\nu\cdot \nabla \Gamma(x^*,\cdot)  d\HH^{n-1}-\frac{4}{3}\int \phi\nu\cdot\nabla \Gamma(x^*,\cdot)  d\HH^{n-1}\\
		=&\phi(x).
	\end{aligned}\]
\end{example}
\subsection{Global lower bounds of the Green function}

The following lemma shows that the even though the indecomposablility assumption in Theorem~\ref{asymp1} is necessary as can be seen from Example \ref{2Planesn-2}, Theorem~\ref{asymp1} may remain to be true for some varifolds with decomposable tangent cones. Examples of minimal submanifolds satisfying Lemma~\ref{regular.at.infinity.2.planes} include the $n$-dimensional catenoid in $\R^{n+1}$ for $n\ge 3$ and the Lawlor's neck introduced in \cite{La}.
\begin{lemma}\label{regular.at.infinity.2.planes}
Let $\Sigma^n\subset \R^{n+k}$ be smooth minimal submanifold with two ends $E_1,E_2$. We assume that for $i=1,2$, $E_i$ is regular in the sense that $E_i$ can be written as a graph of $v_i$ over some $n$-plane $P_i$ in $\R^{n+k}$ where in a coordinate $z$ of $P_i\setminus B_R(0)$ there exists $a,b,c_j\in \R^k$ such that \[v_i(z)=a+b|z|^{2-n}+\sum_{j=1}^nc_jz_j|z|^{-n}+O(|z|^{-n}).\] See \cite[Definition, page 800]{Sch83} for more details. In particular, the tangent cone at infinity of $\Sigma$  is $P_1+ P_2$. We also assume that there exist an isometry of $\Sigma$  such that one end is mapped to the other one.   Then for any  $x\in\Sigma $, \[\lim_{y\to\infty}|x-y|^{n-2}G(x,y)=\frac{1}{2n(n-2)\omega_n}.\] 
\end{lemma}
\begin{proof}
    
Let $u_i=G(x,\cdot)|_{E_i}$, $i=1,2$. Since $\Sigma$ has regular ends, $u_i$ can be transplanted as a function on $P_i\setminus B_{2R}$ and it satisfies a uniformly elliptic equation. By Theorem~\ref{globalGreenvarifolds}, $0\le u_i\le C|x-\cdot |^{2-n}$. Then standard argument as in \cite[Proof of Proposition 3]{Sch83} shows that there exists $c_i\ge 0$ such that
\begin{equation}\label{two.end.example.expansion}
	u_i(z)=c_i|z|^{2-n}+O_1(|z|^{1-n})
\end{equation} where $O_1(|z|^{1-n})$ denotes a function $h$  such that $|h|\le C|z|^{1-n}$ and $|\nabla h|\le C|z|^{-n}$. Since the two ends are symmetric, we have $c_1=c_2=c$. Since $E_i$ are regular, we have $\lim_{|z|\to\infty}|x-y|/|z|= 1$ where $y=(z,v_i(z))$. Thus we get 
\begin{equation}\label{two.end.example.asymptotic}
\lim_{y\to\infty} |x-y|^{n-2}G(x,y)=c.
\end{equation} To get the exact value of $c$, we consider  $b(y)=G(x,y)^{1/(2-n)}$ and the integral \[I(s)=s^{1-n}\int _{\{b=s\}}|\nabla^\Sigma b|=\frac{1}{n-2}\int_{\{b=s\}}|\nabla^\Sigma G|.\] Since $\Sigma$ is a smooth manifold, taking an $s$ derivative, integrating by parts and using $G$ is harmonic away from $x$ (see \cite[(2.6),(2.12)]{CM2} for more details), we have $I(s)$ is a constant.

Since $x\in \Sigma$ is a regular point and $\Theta(\Sigma,x)=1$,  by Theorem~\ref{asymp1} and elliptic regularity \cite{GT}, we have $b(y)\sim (\frac{1}{n(n-2)\omega_n})^{1/(2-n)}|x-y|$ and $|\nabla^\Sigma G(y)|\sim \frac{1}{n\omega_n}|x-y|^{1-n}$ as $y\to x$. Thus
\[\lim_{s\to0}I(s)=\frac{1}{n-2}.\]
By (\ref{two.end.example.expansion}) and (\ref{two.end.example.asymptotic}), we have $b\sim c^{1/(2-n)}|x-y|$ and $|\nabla^\Sigma G|\sim (n-2)c|x-y|^{1-n}$ as $y\to\infty$. Thus 
\[\lim_{s\to\infty}I(s)=\frac{1}{n-2}(\int _{\{b=s\}\cap E_1}|\nabla^\Sigma G|+\int_{\{b=s\}\cap E_2}|\nabla^\Sigma G|)=2cn\omega_n.\]
Thus we get $c=\frac{1}{2n(n-2)\omega_n}$. This finishes the claim.
\end{proof}
	



	\section*{Acknowledgement}

The author would like to express his sincere gratitude to his advisor, Richard Schoen, for constant support and many enlightening discussions. He would also like to thank Chao Li for many helpful  discussions and advice.



\begin{thebibliography}{99}
	
	\bibitem{allard1972first}
	Allard, W. K.: On the first variation of a varifold.  
	Annals of Mathematics 95(3), 417-491 (1972)
	
	
	\bibitem{bombieri1972harnack}
	Bombieri, E., Giusti, E.: Harnack's inequality for elliptic differential equations on minimal surfaces.  
	Inventiones Mathematicae 15(1), 24-46 (1972)
	
	\bibitem{burch1978dini}
	Burch, C.: The Dini condition and regularity of weak solutions of elliptic equations.  
	Journal of Differential Equations 30(3), 308-323 (1978)
	\bibitem{CLY}
	Cheng, S.-Y., Li, P., Yau, S.-T.: Heat equations on minimal submanifolds and their applications.  
	American Journal of Mathematics 106(5), 1033-1065 (1984)
	
	\bibitem{CM1}
	Colding, T. H., Minicozzi, W. P. II: Large scale behavior of kernels of Schr\"odinger operators.  
	American Journal of Mathematics 119(6), 1355-1398 (1997)
	
	\bibitem{CM2}
	Colding, T. H.,  Minicozzi W. P. II: Harmonic functions with polynomial growth. Journal of Differential Geometry 46(1) : 1-77 (1997).
	
	\bibitem{DeLellis2016}
	De Lellis, C.: The size of the singular set of area-minimizing currents.  
	Surveys in Differential Geometry 21(1), 1-83 (2016)
	
	\bibitem{D}
	Ding, Y.: Heat kernels and Green's functions on limit spaces.  
	Communications in Analysis and Geometry 10(3), 475-514 (2002)
	
	\bibitem{gariepy2015measure}
	Evans, L. C., Gariepy, R. F.: Measure theory and fine properties of functions, revised edition.  
	Studies in Advanced Mathematics, CRC Press, Boca Raton, FL (2015)
	
	\bibitem{Fed}
	Federer, H.: Geometric measure theory.  
	Springer (2014)
	
	\bibitem{Fed70}
	Federer, H.: The singular sets of area minimizing rectifiable currents with codimension one and of area minimizing flat chains modulo two with arbitrary codimension.  
	Annals of Mathematics 92, 767-771 (1970)
	
	\bibitem{GT}
	Gilbarg, D., Trudinger, N.: Elliptic partial differential equations of second order.  
	Springer, Vol. 224, No. 2 (1977)
	
	\bibitem{gruter1982green}
	Gr\"uter, M., Widman, K.: The Green function for uniformly elliptic equations.  
	Manuscripta Mathematica 37(3), 303-342 (1982)
	
	\bibitem{La}
	Lawlor, G.: The angle criterion.  
	Inventiones Mathematicae 95(2), 437-446 (1989)

    
	\bibitem{littman1963regular}
	Littman, W., Stampacchia, G., Weinberger, H.: Regular points for elliptic equations with discontinuous coefficients.  
	Annali della Scuola Normale Superiore di Pisa - Classe di Scienze 17(1-2), 43-77 (1963)
	
	\bibitem{LT}
	Li, P., Tam, L.-F.: Complete surfaces with finite total curvature.  
	Journal of Differential Geometry 33(1), 139-168 (1991)
	
	\bibitem{LY}
	Li, P., Yau, S.-T.: On the parabolic kernel of the Schr\"odinger operator.  
	Acta Mathematica 156(3-4), 153-201 (1986)
	
    \bibitem{menne09isoperi}Menne, U.: Some applications of the isoperimetric inequality for integral varifolds. Adv. Calc. Var. 2, 247-269 (2009)
	
    \bibitem{menne2012decay}
	Menne, U.: Decay estimates for the quadratic tilt-excess of integral varifolds.  
	Archive for Rational Mechanics and Analysis 204, 1-83 (2012)
	
	\bibitem{Me1sobolev}
	Menne, U.: Sobolev functions on varifolds.  
	Proceedings of the London Mathematical Society 113(6), 725-774 (2016)
	
	\bibitem{Me2weaklydifferentiable}
	Menne, U.: Weakly differentiable functions on varifolds.  
	Indiana University Mathematics Journal 65(3), 977-1088 (2016)
	
	\bibitem{MS}
	Michael, J. H., Simon, L.: Sobolev and mean value inequalities on generalized submanifolds of $\mathbb{R}^n$.  
	Communications on Pure and Applied Mathematics 26(3), 361-379 (1973)
	
    \bibitem{Modino}Mondino A.: Existence of integral-varifolds minimizing and in Riemannian manifolds. Calculus of Variations and Partial Differential Equations 49(1), 431-70 (2014).
	
    \bibitem{morgan2016gmt}
	Morgan, F.: Geometric measure theory: A beginner's guide.  
	Academic Press, Amsterdam, 5th edition (2016)
	
	\bibitem{petersen2006riemannian}
	Petersen, P.: Riemannian geometry.  
	Springer, Volume 171 (2006)
    
	\bibitem{pitts}Pitts, J. T.: Existence and regularity of minimal surfaces on Riemannian manifolds, Bulletin of the American Mathematical Society 82, no. 3 : 503-504 (1976).
    
    \bibitem{Protter60}
    Protter, M. H.: Unique continuation for elliptic equations, Trans. Amer. Math.
Soc. 95, 81–91 (1960). 
	
    \bibitem{Sch83}
	Schoen, R.: Uniqueness, symmetry, and embeddedness of minimal surfaces. Journal of Differential Geometry 18(4) : 791-809 (1983).
	
	\bibitem{simon1983survey}
	Simon, L.: Survey lectures on minimal submanifolds. Ann. Math. Stud. \textbf{103}, 3-52 (1983)
	
	\bibitem{simon1987asymptotic}
	Simon, L.: Asymptotic behaviour of minimal graphs over exterior domains. Ann. Henri Poincar\'e \textbf{4(3)}, 231-242 (1987)
	
	\bibitem{Sim}
	Simon, L.: Lectures on geometric measure theory.  
	The Australian National University, Mathematical Sciences Institute, Centre for Mathematics \& its Applications (1983)
	
	\bibitem{simon2014introduction}
	Simon, L.: Introduction to geometric measure theory.  
	Tsinghua Lectures (2014)
	
	\bibitem{Var}
	Varopoulos, N. T.: Green's functions on positively curved manifolds.  
	Journal of Functional Analysis 45(1), 109-118 (1982)
	
	\bibitem{Wang}
	Wang, Z.: Mean convex smoothing of mean convex cones. Geometric and Functional Analysis 34(1), 263-301 (2024).
	
	\bibitem{Yau}
	Yau, S.-T. (ed.): Seminar on differential geometry. No. 102. Princeton University Press (1982)
\end{thebibliography}
\end{document}